\documentclass[a4paper,11pt,reqno]{article}
\usepackage{latexsym, amsmath, amsfonts, amssymb, amsthm, amscd, epsfig}
\usepackage{url}

\usepackage[font=small, labelfont={sc}]{caption}

\usepackage[a4paper,scale={0.73,0.75},marginratio={1:1},footskip=10mm,headsep=10mm]{geometry}

\usepackage{color}

\usepackage{hyperref}

\numberwithin{equation}{section}

\def \beq {\begin{eqnarray}}
\def \eeq {\end{eqnarray}}
\def \beqn {\begin{eqnarray*}}
\def \eeqn {\end{eqnarray*}}

\newtheorem{theorem}{Theorem}[section]
\newtheorem{itlemma}[theorem]{Lemma}
\newtheorem{itproposition}[theorem]{Proposition}
\newtheorem{itcorollary}[theorem]{Corollary}
\newtheorem{itremark}[theorem]{Remark}
\newtheorem{itdefinition}[theorem]{Definition}
\newtheorem{itexample}[theorem]{Example}
\newtheorem{itclaim}[theorem]{Claim}
\newtheorem{itfact}[theorem]{Fact}
\newtheorem{itassumption}[theorem]{Assumption}

\newenvironment{fact}{\begin{itfact}\rm}{\end{itfact}}
\newenvironment{claim}{\begin{itclaim}\rm}{\end{itclaim}}
\newenvironment{lemma}{\begin{itlemma}}{\end{itlemma}}
\newenvironment{remark}{\begin{itremark}\rm}{\end{itremark}}
\newenvironment{corollary}{\begin{itcorollary}}{\end{itcorollary}}
\newenvironment{proposition}{\begin{itproposition}}{\end{itproposition}}
\newenvironment{definition}{\begin{itdefinition}\rm}{\end{itdefinition}}
\newenvironment{example}{\begin{itexample}\rm}{\end{itexample}}
\newenvironment{assumption}{\begin{itassumption}}{\end{itassumption}}

\newcommand{\be}[1]{\begin{equation}\label{#1}}
\newcommand{\ee}{\end{equation}}
\newcommand{\bl}[1]{\begin{lemma}\label{#1}}
\newcommand{\br}[1]{\begin{remark}\label{#1}}
\newcommand{\brs}[1]{\begin{remarks}\label{#1}}
\newcommand{\bt}[1]{\begin{theorem}\label{#1}}
\newcommand{\bd}[1]{\begin{definition}\label{#1}}
\newcommand{\bp}[1]{\begin{proposition}\label{#1}}
\newcommand{\bc}[1]{\begin{corollary}\label{#1}}
\newcommand{\bfact}[1]{\begin{fact}\label{#1}.}
\newcommand{\bex}[1]{\begin{example}\label{#1}}
\newcommand{\ec}{\end{corollary}}
\newcommand{\efact}{\end{fact}}
\newcommand{\eex}{\end{example}}
\newcommand{\el}{\end{lemma}}
\newcommand{\er}{\end{remark}}
\newcommand{\ers}{\end{remarks}}
\newcommand{\et}{\end{theorem}}
\newcommand{\ed}{\end{definition}}
\newcommand{\ep}{\end{proposition}}
\newcommand{\epr}{\end{proof}}
\newcommand{\bpr}{\begin{proof}}
\newcommand{\bcl}[1]{\begin{claim}\label{#1}}
\newcommand{\ecl}{\end{claim}}
\newcommand{\bas}[1]{\begin{assumption}\label{#1}}
\newcommand{\eas}{\end{assumption}}

\newcommand{\ecs}{\end{corollary}}
\newcommand{\eers}{\end{exercise}}
\newcommand{\eexs}{\end{example}}
\newcommand{\eems}{\end{example}}
\newcommand{\els}{\end{lemma}}
\newcommand{\eles}{\end{lemmaex}}
\newcommand{\ets}{\end{theorem}}
\newcommand{\eds}{\end{definition}}
\newcommand{\eps}{\end{proposition}}

\newcommand{\bi}{\begin{itemize}}
\newcommand{\ei}{\end{itemize}}
\newcommand{\ben}{\begin{enumerate}}
\newcommand{\een}{\end{enumerate}}

\def\vbar{\mathchoice{\vrule height6.3ptdepth-.5ptwidth.8pt\kern-.8pt}
   {\vrule height6.3ptdepth-.5ptwidth.8pt\kern-.8pt}
   {\vrule height4.1ptdepth-.35ptwidth.6pt\kern-.6pt}
   {\vrule height3.1ptdepth-.25ptwidth.5pt\kern-.5pt}}
\def\fudge{\mathchoice{}{}{\mkern.5mu}{\mkern.8mu}}
\def\bbc#1#2{{\rm \mkern#2mu\vbar\mkern-#2mu#1}}
\def\bbb#1{{\rm I\mkern-3.5mu #1}}
\def\bba#1#2{{\rm #1\mkern-#2mu\fudge #1}}
\def\bb#1{{\count4=`#1 \advance\count4by-64 \ifcase\count4\or\bba A{11.5}\or
   \bbb B\or\bbc C{5}\or\bbb D\or\bbb E\or\bbb F \or\bbc G{5}\or\bbb H\or
   \bbb I\or\bbc J{3}\or\bbb K\or\bbb L \or\bbb M\or\bbb N\or\bbc O{5} \or
   \bbb P\or\bbc Q{5}\or\bbb R\or\bbc S{4.2}\or\bba T{10.5}\or\bbc U{5}\or
   \bba V{12}\or\bba W{16.5}\or\bba X{11}\or\bba Y{11.7}\or\bba Z{7.5}\fi}}

\makeatletter
\newtheorem*{rep@theorem}{\rep@title} \newcommand{\newreptheorem}[2]{%
\newenvironment{rep#1}[1]{%
\def\rep@title{\bf #2 \ref{##1}}%
\begin{rep@theorem} }%
{\end{rep@theorem} } }
\makeatother
\newreptheorem{theorem}{Theorem}
\newreptheorem{lemma}{Lemma}

\def \R {{\mathbb R}}

\def \C {{\cal{C}}}
\def \Q {{\cal{Q}}}
\def \N {{\mathbb N}}
\def \PR {{\mathbb P}}
\def \E {{\mathbb E}}

\def \DD {{\mathbb D}}
\def \s {y}

\def \A {{\cal{A}}}

\newcommand{\ba}[1]{\addtocounter{for}{1} \begin{eqnarray}\label{#1}}
\newcommand{\ea}{\end{eqnarray}}

\def\sqr#1#2{{\vcenter{\vbox{\hrule height .#2pt
                             \hbox{\vrule width .#2pt height#1pt \kern#1pt
                                   \vrule width .#2pt}
                             \hrule height .#2pt}}}}

\def\pmb#1{\setbox0=\hbox{#1}%
   \kern-.025em\copy0\kern-\wd0
   \kern.05em\copy0\kern-\wd0
   \kern-.025em\raise.0433em\box0 }
\def\sqr#1#2{{\vcenter{\vbox{\hrule height.#2pt
     \hbox{\vrule width.#2pt height#1pt \kern#1pt
   \vrule width.#2pt}\hrule height.#2pt}}}}

\def\ve{\varepsilon}
\def\s{\sigma}
\def\d{\delta}
\def\l{\lambda}

\def\g{\gamma}
\def\G{\Gamma}
\def\a{\alpha}
\def\th{\theta}
\def\b{\beta}
\def\r{\rho}

\def\Th{\Theta}

\def\cal{\mathcal}

\newenvironment{myenumerate}{%
\begin{list}{\labelenumi}
	{%
	\setlength{\itemsep}{0.4em}%
	\setlength{\topsep}{0.5em}%
	\setlength\leftmargin{2.6em}%
	\setlength\labelwidth{2.15em}%
	\setlength{\labelsep}{0.45em}%
	\usecounter{enumi}%
	}%
	}%
{\end{list}
}

\renewenvironment{enumerate}{
\renewcommand{\theenumi}{\arabic{enumi}}%
\renewcommand{\labelenumi}{{\rm(\theenumi)}}%
\begin{myenumerate}}%
{\end{myenumerate}}

{\end{myenumerate}}

{\end{myenumerate}}

\newenvironment{myitemize}{%
\begin{list}{$\bullet$}%
 	{%
	\setlength{\itemsep}{0.4em}%
	\setlength{\topsep}{0.5em}%
	\setlength\leftmargin{2.6em}%
	\setlength\labelwidth{2.15em}%
	\setlength{\labelsep}{0.45em}%
	}%
	}%
{\end{list}}

\renewenvironment{itemize}{
\begin{myitemize}}%
{\end{myitemize}}

\newcommand{\Keywords}[1]{\par\noindent 
{\small{\em Keywords\/}: #1}}

\title{Tube estimates for diffusion processes under a weak H\"{o}rmander condition}  

\author{
\textsc{Paolo Pigato}\thanks{%
INRIA, Villers-l\`es-Nancy, F-54600, France.
Universit\'e de Lorraine, IECL, UMR 7502, Vandoeuvre-l\`es-Nancy, F-54600, France. Email: \texttt{%
paolo.pigato@inria.fr}. }
}

\begin{document}

\date{\today}

\maketitle

\begin{abstract} \noindent
We consider a diffusion process under a local weak H\"{o}rmander condition on the coefficients. We find Gaussian estimates for the density in short time and exponential lower and upper bounds for the probability that the diffusion remains in a small tube around a deterministic trajectory (skeleton path), explicitly depending on the radius of the tube and on the energy of the skeleton path. We use a norm which reflects the non-isotropic structure of the problem, meaning that the diffusion propagates in $\R^2$ with different speeds in the directions $\s$ and $[\s,b]$. We establish a connection between this norm and the standard control distance.
\end{abstract}

\smallskip
\Keywords{Density estimates, tube estimates, hypoellipticity, H\"{o}rmander condition, Malliavin Calculus}

\tableofcontents

\section{Introduction}
In this article we consider the following stochastic differential equation on $[0,T]$:
\be{eqn-intro}
X_t=x_0+\int_0^t\s(X_s)\circ dW_s + \int_0^t b(X_s) ds
\ee
where the diffusion $X$ is two-dimensional and the Brownian Motion $W$ is one-dimensional. $\circ dW_s$ denotes the Stratonovich integral, and we suppose a certain geometric property for the diffusion coefficient (which holds true in particular for the equation associated with the Asian option). 
Since $\s$ is just a column vector, the ellipticity assumption fails at any point, and the strong H\"{o}rmander condition fails as well, so we  investigate the regularity of this process assuming a hypoellipticity condition of weak H\"ormander type.
The prototype of this kind of problems is a two dimensional system where the first component $X^1$ follows a stochastic dynamic, and the second component $X^2$ is a deterministic functional of $X^1$, so the randomness acts indirectly on $X^2$. 
Besides the natural application to the Asian option, there are others such as in \cite{HLT:1}, \cite{HLT:2}. In these papers the functioning of a neuron is modeled: $X^2$ is the concentration of some chemicals resulting from a reaction involving the first component $X^1$. Differently from our setting, though, there are several measurements corresponding to the input $X^1$, so $X^2$ is multi-dimensional. The pattern, however, is similar.

We find Gaussian estimates for the density in short time, supposing the process satisfies a weak H\"{o}rmander condition. 
Ben Arous and L\'eandre investigate the decay of the heat kernel of a hypoelliptic diffusion over the diagonal in their celebrated paper \cite{BenArousLeandreII:91}. 
Their framework is different because they work under a strong H\"{o}rmander condition and because they are interested in asymptotic results, whereas we provide results holding for finite positive times.
In \cite{KusuokaStroock:87} explicit two-sided bounds for the density of diffusion processes are established under strong H\"{o}rmander conditions, if the drift is generated by the vector fields of the diffusive part. 
On the opposite, the problem we consider here is of weak H\"{o}rmander type, meaning that the drift has a key role in the propagation of the noise. In this case, the drift gives an additional specific contribution which is usually difficult to handle when trying to estimate the density of the solution. In \cite{BallyKohatsu:10} and \cite{DelarueMenozzi:10} bounds are provided for the density of the Asian type SDE and for a chain of SDEs, in a weak H\"{o}rmander framework. An analytical approach to a similar density estimate is given by Polidoro, Pascucci and Boscain in \cite{Polidoro:97}, \cite{PascucciPolidoro:06}, \cite{BoscainPolidoro:07}.

In this paper, we obtain a more general result than those known in the cited literature, as we allow for a more general coefficient for the Brownian Motion. Indeed we suppose that locally the vector field $\s$ has the same direction of the directional derivative $\partial_\s \s$, whereas the works mentioned above would apply for $\s=(\s_1,0)$ which is a more restrictive condition. Moreover, our coefficients are just locally hypoelliptic.
The other novelty is that thanks to our short time non-asymptotic result we are able to find exponential lower and upper bounds for the probability that the diffusion remains in a small tube around a deterministic trajectory. More precisely we consider \eqref{eqn-intro} and introduce the associated skeleton path solution of the following ODE: 
\[
x_t(\phi)=x_0+\int_0^t\sigma(x_s(\phi))\phi_s ds + \int_0^t b(x_s(\phi)) ds,
\]
for a control $\phi\in L^2[0,T]$. We assume the following weak H\"{o}rmander condition: $\s,\,[\s,b]$ span $\R^2$ locally around $x(\phi)$. 
This is enough to ensure the existence of the density in the case of diffusions (see \cite{Nualart:06}, \cite{Shigekawa:04}). Similar results are also available for SDEs with coefficients with dependence on time, under very weak regularity assumptions (\cite{CattiauxMesnager:02}), SDEs driven by a fractional Brownian Motion (\cite{BaudoinHairer:07}) and for rough differential equations (\cite{CassFriz:10}).

We prove here a tube estimate for \eqref{eqn-intro}, meaning that we find upper and lower bounds for $\PR\left(\sup_{0\leq t\leq T} \|X_t-x_t(\phi)\|\leq R \right)$, explicitly depending on the energy of the skeleton path and on the radius of the tube, that can be time-dependent. Several works have considered this subject, starting from Stroock and Varadhan in \cite{StroockVaradhan:72}, where such result is used to prove the support theorem for diffusion processes. In their work $\|\cdot\|$ is the Euclidean norm, but later on different norms have been used to take into account the regularity of the trajectories (about this, see for example \cite{BenArousGradinaruLedoux:1994} and \cite{FrizLyonsStrook:06}). 
This problem is interesting for physicists because of the Onsager-Machlup functional (see \cite{IkedaWatanabe:89}, \cite{capitaine:00}), and is also related to large and moderate deviation theory (see \cite{bismut1984large}, \cite{guillin2003}).

Since we work under H\"ormander-type conditions, in order to give accurate estimates we consider a norm accounting for the non-diffusive time scale of the process. Indeed, thanks to the H\"ormander condition, the noise propagates in the whole $\R^2$, but with with speed $t^{1/2}$ in the direction $\s$ and $t^{3/2}$ in the direction $[\s,b]$. We also introduce a suitable control metric, adapting the classic control-Carath\'eodory distance, which is equivalent to this norm.

We apply techniques based on the recent work by Bally and Caramellino (\cite{BallyCaramellino:11}, \cite{BallyCaramellino:12}, \cite{BC14}) on density estimates for random variables. In Section 3 we recall some of these results and derive an upper and a lower bound for the density in a fairly abstract framework, starting from the Malliavin-Thalmaier representation formula for the density. The importance of these abstract estimates may go beyond our particular problem.

This paper is organized as follows. In Section 2 we introduce notations and state our main results: the short-time density estimate and the tube estimate. In Section 3 we develop the 
Malliavin calculus techniques that we apply to estimate the density  of our diffusion. In Section 4 we apply these techniques, finding the short-time density estimates mentioned above. In Section 5 we use the short-time result and a concatenation procedure to prove the tube estimate.

\section{Notations and results}

\subsection{Notations}\label{notations}

We start introducing some notations. 
We write $\alpha =(\alpha _{1},...,\alpha _{k})\in \{1,...,n\}^{k}$ for a multi-index with length $\left\vert \alpha \right\vert =k$ and $\partial _{x}^{\alpha }=\partial _{x_{\alpha _{1}}}...\partial _{x_{\alpha_{k}}}$.
For $f,g:\R^n\rightarrow \R^n$ we recall the definition of the directional derivative of $f$ in the direction $g$ as
\[
\partial_g f(x)=(\nabla f)\,g(x) =\sum_{i=1}^n g^i(x)\partial_{x_i}f(x).
\]
The \emph{Lie bracket} $[f,g]$ in $x$ is defined as
\[
[f,g](x)=\partial_f g(x)-\partial_g f(x).
\]
We denote by $M^T$ the transpose of a $2\times 2$ matrix $M$. We also use the notation $\l_*(M)$ for the smallest singular value of $M$, and $\l^*(M)$ for the largest one. We recall that singular values are the square roots of the eigenvalues of $M M^T$, and that, when $M$ is symmetric and semi-definite, singular values coincide with the eigenvalues of $M$. In particular, when $M$ is a covariance matrix, $\l_*(M)$ and $\l^*(M))$ are the smallest and the largest eigenvalues of $M$.

If $M$ is invertible we also associate to $M$ the norm on $\R^2$
\[
|\xi|_M =\sqrt{\langle(M M^T)^{-1} \xi,\xi \rangle}=|M^{-1} \xi|
\] 
For two $2\times 2$ positive semi-definite symmetric matrices $B_1,B_2$, we write $B_1 \leq B_2$ for 
\[
\xi^T B_1 \xi \leq \xi^T B_2 \xi, \quad\mbox{  for all } \xi \in \R^2.
\] 
As we said, we consider the diffusion
\be{eqn}
X_t=x_0+\int_0^t\s(X_s)\circ dW_s + \int_0^t b(X_s) ds,
\ee
where $X$ is in dimension two, $W$ is in dimension one. 
For $x\in \R^2$, we set 
\be{defA}
A(x)=\left(\s(x),[\s,b](x)\right)
\ee
and, for any $R>0$,
\be{defAR}
A_R(x)= \left(R^{1/2}\sigma(x),R^{3/2}[\s,b](x)\right)
\ee

\subsection{Density estimate}
In the first part of the paper we prove an estimate for the density of the solution of \eqref{eqn}. We consider the following assumptions on the coefficients:
\ben
\item[{\bf A1}] The ``first order'' weak H\"ormander condition holds at the initial point of the diffusion:
\[
\l_*(A(x_0))>0
\]
\item[{\bf A2}] $\s,b\in \C^5(\R^2)$ and there exists a constant $\r>0$ such that, $\forall x \in \R^2$:
\[
\sum_{1\leq |\a| \leq 5}
|\partial_x^{\a} \s(x)|+|\partial_x^{\a} b(x)|
\leq \r
\]
\item[{\bf A3}] There exist a neighborhood $V\subset \R^2$ of $x_0$ and a differentiable scalar function $\kappa_\s:V\rightarrow \R$ such that for all $x\in V$
\be{lambdasigma}
\partial_\s \s(x)=\kappa_\s(x) \s(x).
\ee
We suppose that $\sum_{0\leq |\a| \leq 1} |\partial_x^{\a} \kappa_\s(x_0)| \leq \r $.
If $\s(x)=(\s_1(x),0)$, the Asian option stochastic differential equation, this property holds true with $\kappa_\s=\partial_{x_1} \s_1$.  
\een
We prove the following Gaussian bound:
\begin{reptheorem}{mtstime}
Suppose {\bf A1}, {\bf A2}, {\bf A3} hold. Let $(X_t)_{t\in [0,T]}$ be the solution of $\eqref{eqn}$, and for $t\in [0,T]$, let $p_t(x_0,y)$ be the density of $X_t$ at $y$. Then there exist constants $L, C,\d^*$ such that, for any $r>0$, if $0<\d\leq \d^* \exp\left(-L r^2\right)$, setting $\hat{x}_0=x_0+b(x_0)\d$, for $|y-\hat{x}_0|_{A_\d(x_0)} \leq r$
\be{stintro}
\frac{1}{ C \d^2}
\exp\left(-C |y-\hat{x}_0|_{A_\d(x_0)}^2\right)\leq p_\d(x_0,y)
\leq 
\frac{C}{\d^2 }\exp\left(-C^{-1} |y-\hat{x}_0|_{A_\d(x_0)}^2\right)
\ee
\end{reptheorem}
This estimate is local around the point $\hat{x}_0=x_0+\d b(x_0)$. Since we assume the weak H\"ormander condition only at $x_0$, it is not possible to obtain global lower bounds. Indeed the ``local'' weak H\"ormander condition ensures the existence of the density (\cite{KusuokaStroock:84}), but not its positivity. See Example \ref{controllability} for more details on this aspect.

\subsection{Tube estimate}

We suppose $\s,\,b\in \C^5(\R^2)$. For $x\in \R^2$ define 
\[
n(x)=\sum_{k=0}^5 \sum_{|\a|= k} |\partial^\a_x b(x)| + |\partial^\alpha_x \s(x)|,
\]
and set $\l(x)=\l_*(A(x))$. We take now a control $\phi\in L^2[0,T]$, and the associated \emph{skeleton path} solution of 
\be{control}
x_t(\phi)=x_0+\int_0^t\s(x_s(\phi))\phi_s ds + \int_0^t b(x_s(\phi)) ds.
\ee
We denote by $L(\mu,h)$ the class of non-negative functions which have the property
\be{growthcontrol}
f(t)\leq \mu f(s) \quad\mbox{ for }  |t-s| \leq h.
\ee
These functions have been used in \cite{BallyKohatsu:10}, in the choice of an``elliptic evolution sequence",  and in \cite{[BFM]}. They allow us to control the variation of the quantities we are concerned with, along the skeleton path.
In section \ref{diffusion}, when considering the tube estimate, we assume that:
\ben
\item[{\bf H1}]
There exists a function $\l_\cdot:[0,T]\rightarrow(0,1]$ such that
\[
\l(y)\geq \l_t,\quad \forall |y-x_t(\phi)|<1, \quad \forall t\in[0,T].
\]
\item[{\bf H2}] 
There exists a function $n_\cdot:[0,T]\rightarrow[1,\infty)$ such that
\[
n(y)\leq n_t,\quad \forall |y-x_t(\phi)|<1, \quad \forall t\in[0,T].
\]
\item[{\bf H3}] There exists a differentiable scalar function $\kappa_\s:\R^2\rightarrow \R$ s. t.
\[
\partial_\s \s(y)=\kappa_\s(y) \s(y),\quad \forall |y-x_t(\phi)|<1, \quad \forall t\in[0,T]
\]
We suppose also that $|\kappa_\s(y)|\leq n(y),\,|\nabla\kappa_\s (y)|\leq n(y)$.
\item[{\bf H4}] We suppose $|\phi_\cdot|^2,\,\l_\cdot,\,n_\cdot, R_\cdot\in L(\mu,h)$, for some $h>0,\, \mu\geq 1$. 
\een

Notice that the above hypothesis do not involve global controls of our bounds on $\R^2$: they concern the behavior of the coefficients only along the tube, and may vary with $t\in[0,T]$. We stress that also $R_\cdot$, the radius of the tube, may vary with $t$, but that ${\bf H4}$ implies that $\inf_{t\in [0,T]} R_t>0$. This means that we cannot ``squeeze" the tube to $0$ at any time.

For $K, q, K_*, q_* > 0$,  for $0\leq t\leq T$, we denote 
\[
\begin{split}
H_t&= K \left(\frac{\mu n_t}{\l_t}\right)^q,\\
R_t^*(\phi)&=\exp\left(-K_* \left(\frac{\mu n_t}{\l_t}\right)^{q_*}\mu^{2q_*}  \right)
\left( h\wedge \inf_{0\leq \d\leq h} 
\left\{
\d \big/ \int_t^{t+\d} |\phi_s|^2 ds \right\} 
\right).
\end{split}
\]
\begin{reptheorem}{mttubes}
Let $X_t$ be given by \eqref{eqn}, $x_t(\phi)$ by \eqref{control}, and suppose {\bf H1}, {\bf H2}, {\bf H3}, {\bf H4}.
There exist positive constants $K, q, K_*, q_* $ such that, for $H_t$ and $R_t^*(\phi)$ as above, if $R_t \leq R_t ^*(\phi)$  for $0\leq t\leq T$,
\be{tuberesult}
\begin{split}
\exp\left(- \int_0^T H_t
\left(\frac{1}{R_t}+|\phi_t|^2\right)dt \right) &\leq
\PR\left(\sup_{t\leq T} |X_t-x_t(\phi)|_{A_{R_t}(x_t(\phi))}\leq 1 \right) \\
&\quad\quad\quad\quad\quad 
\leq \exp\left(- \int_0^T e^{-H_t}
\left(\frac{1}{R_t}+|\phi_t|^2\right)dt \right).
\end{split}
\ee
In general, even if $R_\cdot$ does not satisfy  $R_t \leq R_t ^*(\phi)$ for $0\leq t\leq T$,  the lower bound holds in the form
\[
\exp\left(- \int_0^T H_t \left( \frac{1}{h}+\frac{1}{R_t}+|\phi_t|^2 dt\right) \right)\leq
\PR\left(\sup_{t\leq T} |X_t-x_t(\phi)|_{A_{R_t}(x_t(\phi))}\leq 1 \right).
\]
\end{reptheorem}
\br{gg}
Notice that estimate \eqref{tuberesult} holds for the controls $\phi$ which belong to the class $L(\mu,h)$, and $\mu$ is involved in the definition of $H_t$. In this sense, $H_t$ depends on the ``growth property" \eqref{growthcontrol} of $\phi$.
\er
Both these theorems can also be stated in a variant of the Carath\'eodory distance which looks appropriate to our framework. Here we just briefly give the definition, for more details see Appendix \ref{controlmetric}.  For $\phi=(\phi^1_s,\phi^2_s)\in L^2((0,1),\R^2)$, set
\[
\|\phi\|_{1,3}^2=
\int_0^1 |\phi_s^1|^2 ds + \left( \int_0^1 |\phi_s^2|^2 ds \right)^{\frac{1}{3}}
\]
and define the class of controls
\[
C_A(x,y)=\{\phi\in L^2((0,1),\R^2): 
dv_s= A (v_s) \phi_s ds,\,x=v_0,\,y=v_1
\}
\]
(recall $A=(\s,[\s,b])$). We set $
d_c(x,y)=\inf 
\left\{ 
\|\phi\|_{1,3}:\phi \in C_A(x,y)\right\}
$.
Just remark that $\|\phi\|_{1,3}$ accounts of the different speed in the $[\s,b]$ direction. We define also the following quasi-distance on $\Omega=\{x\in \R^2:\l_*(A(x))>0\}$. For $x,y\in \Omega$,
\[
d(x,y)<\sqrt{R} \Leftrightarrow |x-y|_{A_R(x)}< 1.
\]
In Appendix \ref{controlmetric} we prove that $d$ and $d_c$ are equivalent quasi-distances, and that Theorem \ref{mttubes} also holds in the following form:
\bc{cor:cartu}
Let $X_t$ be given by \eqref{eqn}, $x_t(\phi)$ by \eqref{control}, and suppose {\bf H1}, {\bf H2}, {\bf H3}, {\bf H4}. There exist constants $C_T>0$ and $R_*>0$ depending on $\s,\,b,\,\mu,\,h$ such that, if $R_t\leq R_*$ for every $t\in [0,T]$, it holds
\begin{multline*}
\exp\left(- C_T \int_0^T \left(\frac{1}{R_t}+|\phi_t|^2\right)dt \right)\leq
\PR\left( d_c (X_t,x_t(\phi))\leq \sqrt{R_t}, \quad \forall t\in [0,T]\right)\\ 
\leq \exp\left(- \frac{1}{C_T}\int_0^T \left(\frac{1}{R_t}+|\phi_t|^2\right)dt \right)
\end{multline*}
\ec

\subsection{Examples and comments}

\bex{controllability}
As mentioned before, assuming the weak H\"{o}rmander condition only in the initial point $x_0$ ensures the existence of the density $p_\d(x_0,y)$, but not its positivity. 
It does not even ensure that the density is positive locally around $x_0$. In \cite{DelarueMenozzi:10}, a multidimensional system under a weak H\"ormander condition is studied, and a global lower bound for the density is provided, but the coefficients are hypoelliptic uniformly on the whole space where the diffusion propagates. 

The fact that we have lower bounds for the density supposing only \textbf{A1} might appear contradictory. In fact, our estimates are local around $\hat{x}_0$, the translated initial condition, and there is no contradiction, as we see in the following classical example (see for instance (3.2.6) in \cite{coron:libro}). Take
\[
X^1_t=1+W_t, \quad\quad
X^2_t=\int_0^t b_2(X_s^1) ds,
\]
where 
\[
b_2(\xi)=\xi^2 1_{\{|\xi|\leq 1\}}+\bar{b}(\xi) 1_{\{|\xi|> 1\}}
\]
and $\bar{b}$ is chosen non-negative and such that $\textbf{A2}$ is satisfied. Weak H\"{o}rmander holds at $X_0=x_0=\left(\begin{array}{c}
1\\
0
\end{array}\right),
$
but for any $y=\left(\begin{array}{c}
y^1\\
y^2
\end{array}\right)$ with $y^2< 0$, $p_\d(x_0,y)=0$, $\forall \d>0$. We have
\[
\s(x_0)=\left(\begin{array}{c}
1\\
0
\end{array}\right),
\quad
b(x_0)=\left(\begin{array}{c}
0\\
(x^1_0)^2
\end{array}\right)
=\left(\begin{array}{c}
0\\
1
\end{array}\right),
\quad
[\s,b](x_0)=\left(\begin{array}{c}
0\\
2 x^1_0
\end{array}\right)
=\left(\begin{array}{c}
0\\
2
\end{array}\right)
\]
In fact, for any fixed $r>0$, the set $\{y:|y-\hat{x}_0|_{A_\d(x_0)}\leq r\}$, on which Theorem \ref{mtstime} holds, is included in $\R\times \R^{+}$, the support of $X_\d$. Indeed $y$ satisfies
\[
|y-\hat{x}_0|_{A_\d(x_0)}=\sqrt{\d^{-1}(y^1-1)^2+\frac{1}{4}\d^{-3}(y^2-\d)^2}\leq r
\]
For $y^2<0$, 
\[
|y-\hat{x}_0|_{A_\d(x_0)}\leq r \Rightarrow \frac{1}{2}\d^{-1/2}\leq r 
\Rightarrow \d\geq \frac{1}{4 r^2}\geq \d^*\exp(-2 L r^2)
\]
if $\d^*\leq \frac{1}{4}$, and this is in contrast with condition $\d\leq \d^*\exp(- L r^2)$ of Theorem \ref{mtstime}.
\eexs

\bex{}
Looking at the geometric condition $\partial_\s \s(x)= \kappa_\s(x) \s(x)$ ({see \bf A3} and {\bf H3}) on the coefficients, it is easy to see that it holds if $\s=(\s_1,0)$. We give here some other simple examples of diffusion coefficient $\s$ satisfying this condition, but with $\s_2\neq 0$:
\bi
\item If $\s=(\s_1,\s_2)$, with $\s_2=C \s_1$ for some constant $C$, we have that the condition is satisfied with $\kappa_\s=\partial_{x_1}\s_1+\partial_{x_2}\s_2$. Remark that with $C=0$ we recover the Asian option SDE.
\item If, for $\a,\b,\g$ constants, 
\[
\s(x_1,x_2)=
\left(\begin{array}{c}
\a x_1 +\b \\
\a x_2 +\g
\end{array}\right)
\]
the condition is satisfied with $\kappa_\s=\a$.
\item If, for $\a, C$ constants,
\[
\s(x_1,x_2)=C
\left(\begin{array}{c}
(x_1/x_2)^\a \\
(x_1/x_2)^{\a-1}
\end{array}\right)
\]
the condition is satisfied with $\kappa_\s=0$.
\ei
These examples show that our estimates are applicable to systems where the regimes of propagation are not completely separated, meaning that the one-dimensional Brownian Motion $W$ can act on both the components of $X$  (improving in this sense the results in \cite{BallyKohatsu:10} and \cite{DelarueMenozzi:10}). On the other hand, the condition required on $\partial_\s \s$ has in some sense the same role of ``separating" the different speeds of propagation. Indeed, we need this assumption to deal with a term of order $t$, which is hard to handle because of its fast speed of propagation, in comparison with the speed $t^{3/2}$ associated to $[\s,b]$.  

For this reason, a multidimensional extension of these results looks quite hard to obtain, especially if we want to consider systems where $W$ is multi-dimensional. This would produce terms of order $t$, associated to the brackets $[\s^i,\s^j]$. To handle these terms we could imagine a generalization of the condition on $\partial_\s \s$, but we believe that this is not an easy task. 
On the other hand, similar results on a multidimensional system, but of strong H\"ormander type, are the subject of the recent work with  Bally and Caramellino (\cite{BCP1,BCP2}), and the techniques used in this paper are also applicable to the system studied in \cite{DelarueMenozzi:10} (cf. \cite{tesipigato}).
\eexs

\bex{}
Consider the geometric Asian option with time horizon $T$ on the  Black \& Scholes model (\cite{FPP}). This can be expressed as
\[
d X^1_t= \s \circ dW_t +r dt=\s dW_t +r dt;\,X^1_0=\xi,\quad\quad\quad 
d X^2_t = \frac{X^1_t}{T} dt;\,X^2_0=0.
\]
In this case, for $R>0$ fixed constant,
\[
A_R^{-1}(x)=
\left(\begin{array}{cc}
\s R^{1/2} & 0 \\
0  & \frac{\s}{T} R^{3/2}      
\end{array}\right)^{-1}
=\frac{1}{\s} \left(\begin{array}{cc}
\frac{1}{R^{1/2}} & 0 \\
0  & \frac{T}{R^{3/2}}      
\end{array}\right)
\]
does not depend on $x$. We take as control
$\phi_t=0$ so $x_t(\phi)= \left(\xi+r t, \frac{\xi t+ r t^2/2}{T}\right)$. We have
\[
\begin{split}
|X_t-x_t(\phi)|_{A_R(x_t(\phi))}
&=
\frac{1}{\s} \sqrt{\frac{|X^1_t-(\xi+r t)|^2}{R}+\frac{T^2 |X^2_t-(\xi t+ r t^2/2)/T|^2}{R^3}}\\
&=
\frac{1}{\s} \sqrt{\frac{|\s W_t|^2}{R}+\frac{|\s\int_0^t W_s ds |^2}{R^3}},
\end{split}
\]
and \eqref{tuberesult} gives
\[
e^{-C_1 T/R}\leq \PR\left( \sup_{t\leq T} \left\{\frac{|W_t|^2}{R}+\frac{|\int_0^t W_s ds |^2}{R^3}\right\} \leq 1 \right)
\leq e^{-C_2 T/R}.
\]
\eexs

\bex{}
Consider a system given  by the Black \& Scholes model for the price of an asset, and an (arithmetic average) Asian option on that asset with time horizon $T$ (see for instance \cite{yor92, carr, FPP}). This is a model of real interest in mathematical finance. The associated SDE is
\[
d X^1_t= X_t^1 (\s \circ dW_t +r dt);\,X^1_0= \xi>0,\quad\quad\quad 
d X^2_t = \frac{X^1_t}{T} dt;\,X^2_0=0,
\]
and $X^1_t=\xi e^{\s W_t+rt}$. The stochastic integral is in Stratonovich form so to recover the classical formulation $r\rightarrow r+\s^2/2$.
In this case, for $R>0$ fixed constant,
\[
A_R^{-1}(x)=
\left(\begin{array}{cc}
\s x^1 R^{1/2} & 0 \\
0  & \frac{\s x^1}{T} R^{3/2}      
\end{array}\right)^{-1}
=\frac{1}{\s x^1} \left(\begin{array}{cc}
\frac{1}{R^{1/2}} & 0  \\
0  & \frac{T}{R^{3/2}}      
\end{array}\right)
\]
Remark that this matrix is invertible for $x^1\neq 0$. Since we are working under local non-degeneracy assumptions, our tube estimates hold for any initial condition $\xi>0$, provided that $R>0$ is small enough, since this implies the positivity of the first component of the skeleton path at any time $t>0$. On the other hand, results requiring ``global" non degeneracy, such as the density estimates in \cite{DelarueMenozzi:10}, do not hold for this model. We take as control
$\phi_t=0$ so $x_t(\phi)=\xi \left( e^{r t} , \frac{1}{T} \int_0^t e^{rs} ds \right)$. We have
\[
\begin{split}
|X_t-x_t(\phi)&|_{A_R(x_t(\phi))}
=
\frac{1}{\s \xi e^{r t} } \sqrt{\frac{|X^1_t-\xi e^{r t}|^2}{R}
+\frac{T^2 |X^2_t- \frac{\xi}{T} \int_0^t e^{rs} ds |^2}{R^3}}\\
&=
\frac{1}{\s \xi e^{r t} } \sqrt{\frac{\xi^2| e^{r t}(e^{\s W_t}-1)|^2}{R}
+\frac{\xi^2 |\int_0^t  e^{r s +\s W_s} ds
- \int_0^t e^{r s} ds|^2}{R^3}}\\
&=
\frac{1}{\s e^{r t}} \sqrt{\frac{|e^{r t}(e^{\s W_t}-1)|^2}{R}
+\frac{|\int_0^t e^{r s}(e^{\s W_s}-1) ds|^2}{R^3}}\\
\end{split}
\]
and \eqref{tuberesult} gives
\[
e^{-C_1 T/R}\leq \PR\left( \sup_{t\leq T} 
\left\{
\frac{|e^{\s W_t}-1|^2}{R \s^2}
+\frac{|\int_0^t e^{r (s-t)}(e^{\s W_s}-1) ds|^2}{R^3 \s ^2} \right\}
\leq 1 \right)
\leq e^{-C_2 T/R}.
\]
\eexs

\section{Malliavin calculus and density estimates}\label{sec:malliavin}

\subsection{Notations}
Our main reference for this section is \cite{Nualart:06}. We consider a probability space $(\Omega,\mathcal{F},\PR)$ and a Brownian motion $W=(W^1_t,...,W^d_t)_{t\geq 0}$.
We denote by $\DD^{k,p}$ the space of the random variables which are $k$ times differentiable in the Malliavin sense in $L^p$, and $\DD^{k,\infty}=\bigcap_{p=1}^\infty \DD^{k,p}$. For a multi-index $\a=(\a_1,\dots, \a_m)$ we denote by $D^\a F$ the Malliavin derivative of $F$ corresponding to the multi-index $\a$.

$\DD^{k,p}$ is the closure of the space of the simple functionals with respect to the
Malliavin Sobolev norm
$$
\|F\|_{k,p}=[\E|F|^p+\sum_{j=1}^k \E |D^{(j)} F|^p]^\frac{1}{p}
$$
where
$$
|D^{(j)} F|=\left(\sum_{|\alpha|=j}\int_{[0,T]^j}|D^\alpha_{s_1,...,s_j}F|^2 d s_1 ... d s_j\right)^{1/2}.
$$
For the special case $j=1$, we use the standard notation
$$
|D F|=|D^{(1)} F|=\left(\sum_{m=1}^d \int_{[0,T]}|D^m_s F|^2 d s \right)^{1/2}.
$$
Hereafter, for $j\in \N \setminus\{0\}$, we write $D^{(j)}$ for the ``derivative of order $j$" and $D^j$ for the ``derivative with respect to $W^j$".

As usual, we also denote by $L$ the Ornstein-Uhlenbeck operator, i.e. $L=-\d\circ D$, where $\d$ is the adjoint operator of $D$.

For a random vector $F=(F_1,...,F_n)$ in the domain of $D$, we define its \emph{Malliavin covariance matrix} as follows:
\[
\g_F^{i,j}=\langle D F_i,D F_j\rangle_\mathcal{H} = \sum_{k=1}^d \int_0^T D^k_s F_i\times D^k_s F_j ds.
\] 
We say that $F$ is \emph{non-degenerate} if its Malliavin covariance matrix is invertible and
\be{nondeg}
\E(|\det \g_F|^{-p})<\infty,\quad \forall p\in \N.
\ee
We denote by $\hat{\g}_F$ the inverse of $\g_F$.

\subsection{Localization}\label{local}
The following notion of localization is introduced in \cite{BallyCaramellino:12}. Consider a random variable $U\in[0,1]$ and denote
$$
d\PR_U=Ud\PR.
$$
$\PR_U$ is a non-negative measure (not a probability measure, in general). We also set $\E_U$ the expectation (integral) w.r.t. $\PR_U$, and denote
\begin{align*}
\| F \|_{p,U}^p&=\E_U(|F|^p)=\E(|F|^p U)\\
\| F \|_{k,p,U}^p&=\| F \|_{p,U}^p+\sum_{j=1}^k \E_U(|D^{(j)}F|^p).
\end{align*}
We assume that $U\in \DD^{2,\infty}$ and for every $p\geq 1$ 
\[
m_U(p):=1+ (\E_U | D \ln U |^{p})^{1/p} + (\E_U | D^{(2)} \ln U |^{p})^{1/p}< \infty.
\]
(notice that our definition of $m_U$ is slightly different from the definition in \cite{BallyCaramellino:12}: we are taking $p$-norms instead of moments, and we also consider $D^{(2)}$, whereas in \cite{BallyCaramellino:12} only the first order derivative $D$ appears in $m_U$).
For $F=(F^1,\cdots,F^n)$ such that $F^1,\cdots,F^n\in \DD^{2,\infty}$ and $V\in \DD^{1,\infty}$, for any localization function $U$ we introduce the localized Malliavin weights
\[
H_{i,U}(F,V)=
\sum_{j=1}^n V \hat{\g}_F^{i,j} L F^j -\langle D(V\hat{\g}_F^{i,j}) , DF^j \rangle
-V \hat{\g}_F^{i,j} \langle D\ln U , DF^j \rangle
\]
and the vector
\[
H_{U}(F,V)=\left(H_{i,U}(F,V)\right)_{i=1,\dots n}.
\]
The following \emph{representation formula for the localized density} has been proved in \cite{BallyCaramellino:11}.
\bt{bhh}
Let $U$ be a localizing r.v. such that under $\PR_U$ \eqref{nondeg} holds, i.e.
\[
\E_U[|\det \g_F|^{-p}]<\infty,\quad \forall p\in \N.
\]
Then, under $\PR_U$ the law of $F$ is absolutely continuous and has a continuous density $p_{F,U}$ which may be represented as
\be{repden}
p_{F,U}(x)=\sum_{i=1}^n \E_U[\partial_i \Q_n(F-x) H_{i,U}(F,1)]
\ee
where $\Q_n$ denotes the Poisson kernel on $\R^n$, i.e. the fundamental solution of the Laplace operator $\Delta \Q_n=\d_0$. This is given by
\[
\Q_1(x)=\max(x,0);\quad \Q_2(x)=\A_2^{-1}\ln |x|;\quad \Q_n(x)=-\A_n^{-1} |x|^{2-n}, \, n>2,
\]
where $\A_n$ is the area of the unit sphere in $\R^n$.
\et
This is a localized version of the formula
\[
p_F(x)=\sum_{i=1}^n \E\left[\partial_i\Q_n(F-x)H_i(F,1)\right]
\]
where the Malliavin weights are given by
\[
H(F,G) = G \hat{\g}_F\times LF - \langle D(\hat{\g}_F G),DF\rangle
\]
for which we refer to \cite{MalliavinThalmaier:06}. We recall the following relation between localized weights, which can be easily checked (a similar formula is proved in \cite{BallyCaramellino:12}). For any $U,V$ localizing r.v.s, $F,G\in \DD^{2,\infty}$ 
\be{changewei}
H_U(F,VG )=V H_{UV}(F,G)
\ee

\bex{} 
The following example of localizing function is taken from \cite{BallyCaramellino:12}. Consider the function depending on a parameter $a>0$:
\[
\psi_a(x)=1_{|x|\leq a}+\exp\left( 1-\frac{a^2}{a^2-(x-a)^2} \right)1_{a<|x|<2a},
\]
which is a smooth version of the indicator function $1_{\{|x|\leq a\}}$.
For $\Theta_i\in \DD^{1,\infty}, \,i=1\dots n$, and $r>0$, we define the localization r. v.
\be{defU}
U_r=\prod_{i=1}^n \psi_{r}(\Theta_i)
\ee
For this choice of $U_r$ we have that for any $p\geq 1$,
\be{normlocalization}
m_{U_r}(p)
\leq  
C_p \left( 1+\frac{\|\Theta\|_{2,p}^2}{r^2} \right)
\ee
and
\be{inequality2}
\|1-U_r\|_{1,p} \leq  C \left( 1+\frac{\|\Th\|_{1,2p}}{r}\right) \sum_{i=1}^n \PR (|\Th_i|\geq r)^{1/2p}.
\ee
The proof of \eqref{normlocalization} follows from inequalities
\be{boundpsi}
\sup_x |(\ln \psi_a)' (x)|^p \psi_a(x)\leq
\frac{4^{p}}{a^{p}} \sup_{t\geq 0}(t^{2 p} e^{1-t})\leq \frac{C_p}{a^p}<\infty
\ee
and
\be{boundpsi2}
\sup_x |(\ln \psi_a)'' (x)|^p \psi_a(x)
\leq
\frac{8^{p}}{a^{2 p}} \sup_{t\geq 0}(t^{3 p} e^{1-t})
+
\frac{2^{p}}{a^{2 p}} \sup_{t\geq 0}(t^{2 p} e^{1-t})
\leq \frac{C_p}{a^{2p}}
<\infty
\ee
Indeed
\[
\begin{split}
U_r |D \ln U_r|^p
&= \prod_{i=1}^n \psi_{r}(\Theta_i) \Big| \sum_{i=1}^n (\ln \psi_r)' (\Th_i) D\Th_i \Big|^p \\
&\leq
\prod_{i=1}^n \psi_{r}(\Theta_i) \Big( \sum_{i=1}^n | (\ln \psi_r)' (\Th_i)|^2 \Big)^{p/2} \Big( \sum_{i=1}^n | D \Th_i |^2 \Big)^{p/2}\\
&\leq
c_p \Big(  \sum_{i=1}^n | (\ln \psi_r)' (\Th_i)|^p \psi_{r}(\Theta_i)  \Big)
|D\Th |^{p}.
\end{split}
\]
Here we apply \eqref{boundpsi}, and find
\be{eqU}
U_r |D \ln U_r|^p
\leq
C_p \frac{|D\Th|^p}{r^p}.
\ee
This implies $(\E_{U_r} |D \ln U_r|^p)^{1/p}\leq C_p \frac{\|\Th\|_{1,p}}{r}$. 
We also have, using \eqref{boundpsi} and \eqref{boundpsi2},
\[
\begin{split}
U_r |D^{(2)}& \ln U_r|^p
= \prod_{i=1}^n \psi_{r}(\Theta_i) \Big| D \left( \sum_{i=1}^n (\ln \psi_r)' (\Th_i) D\Th_i \right) \Big|^p \\
&\leq C_p \prod_{i=1}^n \psi_{r}(\Theta_i) \left[
\Big|  \sum_{i=1}^n (\ln \psi_r)'' (\Th_i) (D\Th_i)^2 \Big|^p +
\Big|  \sum_{i=1}^n (\ln \psi_r)' (\Th_i) D^{(2)}\Th_i \Big|^p \right]\\
&\leq C_p 
\Big(  \sum_{i=1}^n |(\ln \psi_r)'' (\Th_i)|^p \psi_{r}(\Theta_i)\Big)  |D\Th|^{2p}  +
C_p
\Big(  \sum_{i=1}^n |(\ln \psi_r)' (\Th_i)|^p \psi_{r}(\Theta_i)\Big)   |D^{(2)}\Th|^p \\
&\leq 
C_p \Big( \frac{|D\Th|^{2p}}{r^{2p}}  +\frac{ |D^{(2)}\Th|^p}{r^p} \Big)
\end{split}
\]
and so 
\[
(\E_{U_r} |D^{(2)} \ln U_r|^p)^{1/p} \leq C_p \left(\left( \frac{\|\Th\|_{1,p}}{r} \right)^2+\frac{\|\Th\|_{2,p}}{r} \right). 
\]
This proves \eqref{normlocalization}
Moreover, since $D_s U_r=0$ on $\bigcap_i \{|\Th_i| < r\}=\big( \bigcup_i \{|\Th_i|\geq r\}\big)^c$,
\[
D_s (1-U_r)= - 1_{\{\bigcup_i \{|\Th_i|\geq r\}\}} D_s U_r
\]
and from H\"older inequality
\[
\E |D_s (1-U_r)|^p \leq (\E 1_{\{\bigcup_i \{|\Th_i|\geq r\}\}})^{1/2} (\E |D_s U_r|^{2p})^{1/2}
\]
We control the first factor with the tail estimate
\[
 (\E 1_{\{\cup_i \{|\Th_i|\geq r\}\}})^{1/2} \leq C\sum_{i=1}^n \PR (|\Th_i|\geq r)^{1/2},
\]
and
we also have
\[
| D_s U_r |^{2p} \leq U_r |D \ln U_r|^{2p},
\]
and from \eqref{eqU}
\[
(\E | D_s (1-U_r) |^{p})^{1/p} \leq C_p \frac{\|\Th\|_{1,2p}}{r} \sum_{i=1}^n \PR (|\Th_i| \geq r)^{1/2p}.
\]
Moreover
\[
\E |1-U_r|^p\leq \PR(1-U_r>0) \leq \PR(|\Th_i|>r, \, \exists i=1,\dots n)\leq
\sum_{i=1}^n \PR(|\Th_i|>r),
\]
so \eqref{inequality2} is proved.
\eexs
 
\subsection{The distance between two local densities}
We discuss some techniques, based on Malliavin calculus, for estimating the density of a random variable. These ideas are based on the recent work of Bally and Caramellino (\cite{BallyCaramellino:12}, \cite{BC14}). 

In what follows for a given matrix $A$ we consider its Frobenius norm, given as
\[
\|A\|_{Fr}=\sqrt{\sum_{i,j} |A_{i,j}^2|}=\sqrt{Tr(A^T A)}.
\]
We will employ the fact that the Frobenius norm is sub-multiplicative. 
Take a square $d\times d$ matrix $\g$, symmetric and positive definite. Recall that we denote by $\l^*(\g)$ and $\l_*(\g)$ the largest and the smallest singular values of $\g$, which in this case coincide with the largest and smallest eigenvalues. 
From the equivalence between Frobenius and spectral norm we have 
\[
\l^*(\g)\leq \|\g\|_{Fr} \leq \sqrt{d}  \l^*(\g).
\]
Denoting $\hat{\g}=\g^{-1}$, it holds $\l^*(\hat{\g})=1/\l_*(\g)$. So  
\[
\frac{1}{\l_*(\g)} \leq \| \hat{\g} \|_{Fr} \leq \frac{\sqrt{d}}{\l_*(\g)}.
\]
For two time dependent matrices $A_s,B_s$, we have the following ``Cauchy-Schwartz" inequality:  
\[
\|\int A_s B_s ds\|_{Fr}^2 \leq \int \|A_s\|_{Fr}^2 ds \int \|B_s\|_{Fr}^2 ds.
\]
In particular, if $B_s=v_s$ is a vector,
\[
|\int A_s v_s ds|^2 \leq \int \|A_s\|_{Fr}^2 ds \int |v_s|^2 ds.
\]
We fix some notation. Let $W$ be a Brownian Motion in $\R^d$.
For two random variables $F=(F_1,\dots F_n),\,G=(G_1,\dots G_n)$ in $\DD^{3,\infty}$ and a localizing r. v. $U$, we denote
\[
\begin{split}
\Gamma_{F,U}(p)
&=1+\left( \E_U \l_*(\g_F)^{-p}\right)^{1/p}\\
\Gamma_{F,G,U}(p)
&=1+\sup_{0\leq \ve \leq 1} \left( \E_U \l_*(\g_{G+\ve (F-G)})^{-p}\right)^{1/p}\\
n_{F,G,U}(p)
&= 1+\|F\|_{3,p,U}+\|G\|_{3,p,U}+\|LF\|_{1,p,U}+\|LG\|_{1,p,U}\\
\Delta_2 (F,G)
&=|D(F-G)|+|D^{(2)}(F-G)|+|L(F-G)|
\end{split}
\]
We also write $n_{F,U}(p)$ for $n_{F,0,U}(p)$. Moreover, in all the above notations, when $U=1$, i.e. the localization is ``trivial", we omit it in the notation. 
Remark that notations $n_{F,U}$ and $n_{F,G}$, although similar, denote different things.
Since we are differentiating with respect to a Brownian Motion, as a direct consequence of Meyer's inequality (see for instance \cite{Nualart:06}), we have
\[
n_{F,G,U}(p)
\leq 1+C\left(\|F\|_{3,p}+\|G\|_{3,p}\right)
\]
for every $F,G,U$. 

We now give the main result of this section, comparing  the densities of the laws of two random variables under $\PR_U$.
\bt{genden}
Let $U$ be a localizing r.v. with $m_U(32n )<\infty$. Let $F=(F_1,\dots,F_n),\,G=(G_1,\dots,G_n)\in \DD^{3,32 n}$. Suppose $\Gamma_{G,U}(p)<\infty$ 
and $\Gamma_{F,U}(p)<\infty$ for any $p>1$.
Then there exists a constant $C_1$ such that
\[
p_{G,U}(y)
- C_1 
\|\Delta_2(F,G)\|_{32 n,U} \leq p_{F,U}(y) 
\leq p_F(y)
\]
If, in addition, $\Gamma_{F}(32 n)<\infty$, there exists a constant $C_2$ such that
\[
p_F(y)\leq p_{G,U}(y)
+ C_2 (\|\Delta_2(F,G)\|_{32 n,U}+\|1-U\|_{1,14n})
\]
\et
\br{c}
We can take 
\[
\begin{split}
C_1&=C\left[ m_U(32n ) \Gamma_{G,U}(32 n) n_{F,G,U}(32 n)\right]^{24 n^2}\\
C_2&=C\left[m_U(32n ) \Gamma_{F}(32 n) n_{F,G}(32 n)\right]^{24 n^2}
\end{split}
\]
where $C$ is a constant depending only on the dimension $n$.
\er
The lower bound for $p_{F,U}$ is a version of Proposition 2.5. in
\cite{BallyCaramellino:12}, where here we have specified as possible choice for the exponent $p=32n$. Moreover, we find here that in $m_U$ and $n_{F,G,U}$ we need to consider one more order of derivatives with respect to \cite{BallyCaramellino:12}.
Similar estimates can be found also in \cite{BC14}. 

Before proceeding with the proof we need some preliminary results. We start with an estimate for the localized Malliavin weights and for the difference of weights:
\bl{weightnorm}
Let $U$ be a localizing r.v, $V\in\DD^{1,\infty}$, $F=(F_1,\dots,F_n)\in \DD^{3,\infty}$. Suppose $\Gamma_{F,U}(q)<\infty$ for any $q>1$.
For fixed $p\geq1$, $p_i\geq1,\,i=1,\dots, 4$, with
$\frac{1}{p}=\frac{1}{p_1}+\frac{1}{p_2}+\frac{2}{p_3}+\frac{3}{p_4}$,
there exists a constant $C$ depending only on $p$ and the dimension $n$ such that
\be{wei1} 
\|H_U(F,V)\|_{p,U} \leq C \|V\|_{1,p_1}
 m_U(p_2) \Gamma_{F,U}(p_3)^{2} n_{F,U}(p_4)^{3} 
\ee
Moreover if $\frac{1}{p}=\frac{1}{p_1}+\frac{1}{p_2}+\frac{3}{p_3}+\frac{5}{p_4}$ and $V\in\DD^{2,\infty}$,
\be{wei2} 
\|H_U(F,V)\|_{1,p,U} \leq C \|V\|_{2,p_1}
 m_U(p_2) \Gamma_{F,U}(p_3)^{3} n_{F,U}(p_4)^{5}, 
\ee
Let now $G=(G_1,\dots,G_n)\in \DD^{3,\infty}$. If $\Gamma_{F,G,U}(q)<\infty$ for any $q>1$,  for fixed $p_i\geq 1,\,i=1,\dots, 5$ with
$\frac{1}{p}=\frac{1}{p_1}+\frac{1}{p_2}+\frac{3}{p_3}+\frac{4}{p_4}+\frac{1}{p_5}$, it also holds
\be{wei3}
\left\| H_{U}(F,V)-H_{U}(G,V)\right\|_{p,U}\leq C
\|V\|_{1,p_1} m_U(p_2) \Gamma_{F,G,U}(p_3)^{3} n_{F,G,U}(p_4)^{4}
\|\Delta_{2}(F,G)\|_{p_5,U}.
\ee
\el

\bpr
Consider the weight:
\be{wei}
H_U(F,V)=
V[\hat{\g}_F\times L F -\langle D\hat{\g}_F , DF \rangle]
-\langle \hat{\g}_F  (D V+ VD\ln U) , DF \rangle
\ee
Recall that $D^{(k)}$ means ``derivative of order $k$" and $D^k$ means ``derivative with respect to $W^k$". 
We first consider $D \g_F$ and have the following estimate:
\[
\begin{split}
&\sum_{l=1}^d \int \|D^l_s \g_F\|_{Fr}^2 ds\\
&\quad=
\sum_{l=1}^d \int  \left\| \left( \sum_{k=1}^d \int_0^t  D^l_s D^k_u F_i \times D^k_u F_j +  D^k_u F_i \times D^l_s D^k_u F_j du \right)_{i,j} \right\|_{Fr}^2 ds \\
&\quad\leq
4 | D^{(2)}F |^2 | D F |^2
\end{split}
\]
We now consider $D \hat{\g}_F$. From the chain rule and the derivative of the inversion of matrices,
\be{depr}
D^k \hat{\g}_F = - \hat{\g}_F (D^k \g_F) \hat{\g}_F.
\ee
So, applying also the previous estimate
\[
\sum_{k=1}^d \int \|D^k_s \hat{\g}_F\|_{Fr}^2 ds 
\leq
\|\hat{\g}_F\|_{Fr}^4  \sum_{k=1}^d \int  \|D^k_s \g_F\|_{Fr}^2 ds
\leq
4 \|\hat{\g}_F\|_{Fr}^4 |D F |^2 |D^{(2)} F |^2.
\]
From \eqref{wei} we see that
\[
\begin{split}
|H_U&(F,V)|  \leq
|V| \left( \|\hat{\g}_F \|_{Fr} | L F | 
+ \Big( \sum_{k=1}^d 
\int \| D^k\hat{\g}_F \|_{Fr}^2 ds\Big)^{1/2} \,
| DF | \right)\\
&\quad\quad\quad\quad\quad\quad\quad\quad\quad\quad\quad\quad\quad\quad\quad\quad+
\| \hat{\g}_F \|_{Fr} \big(
|D V|+ | V| |D\ln U| \big)  |DF|\\
&\leq C (|V|+|DV| ) (1+ |D\ln U| )(|DF|+|LF|)
\Big(
\|\hat{\g}_F \|_{Fr} 
+ \Big(\sum_{k=1}^d 
\int \| D^k\hat{\g}_F \|_{Fr}^2 ds\Big)^{1/2}
 \Big)\\
&\leq C (|V|+|DV| ) (1+ |D\ln U| )(1+|DF|+|D^{(2)}F|+|LF|)^3
(1+\|\hat{\g}_F \|_{Fr} )^2
\end{split}
\]
Now 
\[ 
\|H_U(F,V)\|_{p,U} \leq C \|V\|_{1,p_1}
 m_U(p_2) \Gamma_{F,U}(p_3)^{2} n_{F,U}(p_4)^{3} ,
\]
for $\frac{1}{p}=\frac{1}{p_1}+\frac{1}{p_2}+\frac{2}{p_3}+\frac{3}{p_4}$,
follows easily applying H\"older and Minkowski inequalities for $L_p$ norms.

The estimate of $\|H_U(F,V)\|_{1,p,U}$ follows using very similar techniques. 
The part giving the ``main" contribution is $D^{(2)}\hat{\g}_F$, for which, iterating \eqref{depr}, it is not difficult to see
\[
| D^{(2)}\hat{\g}_F |
\leq  C (|D F|+\dots+|D^{(3)}F|)^{4} \,\|\hat{\g}_F\|_{Fr}^{3}
\]
This term is also multiplied by $|DF|$, so we have the estimate of the term giving the main contribution. We leave out the similar estimate of the other terms.

When considering the difference $\|H_U(F,V)-H_U(G,V)\|_{p,U}$, we use similar arguments and the following property of norms: $|ab-cd|\leq |a-c| |b|+|c||b-d|$.
As before the main contribution comes from $D (\hat{\g}_F-\hat{\g}_G)$, so we consider this and leave out the estimates of the other terms.
We remark that
\[
\hat{\g}_F-\hat{\g}_G=\hat{\g}_F(\g_G-\g_F)\hat{\g}_G
\]
and differentiate this product, finding
\[
\begin{split}
&|D (\hat{\g}_F-\hat{\g}_G)|
\leq C
(1+\|\hat{\g}_{F}\|_{Fr}\vee \|\hat{\g}_{G}\|_{Fr})^3\\
& \left(1+|D \g_F|\vee|D \g_G|\right) \left( |\g_F-\g_G|+|D (\g_F-\g_G)| \right) 
\end{split}
\]
where
\[
1+|D \g_F|\vee|D \g_G|
\leq
C \left(1+\sum_{i=1}^{2}|D^{(i)} F|\vee|D^{(i)} G|\right)^2
\]
We have
\[
|\g_F-\g_G| \leq C  |D (F-G)| \,| D (F+G)|
\]
and
\[
|D(\g_F-\g_G) | \leq C  \left(|D(F-G)|+|D^{(2)} (F-G)| \right) \left(|D(F+G)|+|D^{(2)} (F+G)| \right)
\]
Multiplying with $|DF|$, and applying H\"older inequality, we prove the statement.
\epr

\bl{lemma1}
Let $U$ be a localizing r.v.,  $F=(F_1,\dots,F_n), G=(G_1,\dots,G_n)\in \DD^{3,\infty}$. If $\Gamma_{F,G,U}(q)<\infty$ for any $q>1$, there exists a constant $C$ depending only on the dimension $n$ such that
\[
|p_{F,U}(y)-p_{G,U}(y)|
\leq C
\left[m_U(32n ) \Gamma_{F,G,U}(32 n) n_{F,G,U}(32 n)\right]^{12 n^2}
\|\Delta_2(F,G)\|_{32 n,U}
\]
\el

\bpr
We write the densities using \eqref{repden}:
\[
\begin{split}
p_{F,U}(y)-p_{G,U}(y)
&=E_U(
\langle \nabla \cal{Q}_n(F-y),H_U(F,1)\rangle
-\langle \nabla \cal{Q}_n(G-y),H_U(G,1)\rangle)\\
&=E_U\langle \nabla \cal{Q}_n(F-y),H_U(G,1)-H_U(F,1)\rangle
\\
&+E_U \langle \nabla \cal{Q}_n(G-y)-\nabla\cal{Q}_n(F-y),H_U(G,1)\rangle\\
&=
I+J
\end{split}
\]
We recall the following inequality proved in \cite{BallyCaramellino:11}. For $p>n$, with $p'=p/(p-1)$, 
\[
(\E_U|\nabla \cal{Q}_n(F-y)|^{p'})^{1/p'}
\leq C_{p,n} (\E_U |H_U(F,1)|^p )^{p\frac{n-1}{p-n}}.
\]
In particular, for $p=2n$ (fixed from now on), applying \eqref{wei1} with $k=0,p_1=p_2=p_3=p_4=7p=14n$,
\be{estQ}
\begin{split}
&(\E_U|\nabla \cal{Q}_n(F-y)|^{2n/(2n-1)})^{(2n-1)/(2n)}\\
&\quad\leq C (\E_U |H_U(F,1)|^{2n} )^{2(n-1)}\\
&\quad\leq C \left[
m_U(14n) \Gamma_{F,U}(14n)^2 n_{F,U}(14n)^3
\right]^{4n(n-1)}.
\end{split}
\ee
We use now Lemma \ref{weightnorm} to estimate $I$ and $J$. From H\"older inequality
\[
\begin{split}
I=&\E_U |\langle \nabla \cal{Q}_n(F-y),H_U(G,1)-H_U(F,1)\rangle| \\
&\leq
\|\nabla \cal{Q}_n(F-y)\|_{\frac{2n}{2n-1},U}
\|H_U(G,1)-H_U(F,1)\|_{2n,U} 
\end{split}
\]
and we have just provided the estimate for the first factor. For the second we apply \eqref{wei3} with $p_1=p_2=p_3=p_4=p_5=20n$
\[
\begin{split}
&\|H_U(F,1)-H_U(G,1)\|_{2n,U} \\
&\quad\leq C m_U(20n) 
\Gamma_{F,G,U}(20n)^3 n_{F,G,U}(20n)^4
\|\Delta_2(F,G)\|_{20n,U},
\end{split}
\]
We now study $J$. For $\l\in[0,1]$ we denote $F_\l=G+\l(F-G)$. With a Taylor expansion, applying H\"older inequality, integrating again by parts and denoting $V_{j,k}=H_{j,U}(G,1) (F-G)_k$. 
\[
\begin{split}
&\E_U \langle \nabla \cal{Q}_n(F-y)-\nabla \cal{Q}_n(G-y),H_U(G,1)\rangle\\
& \quad =
\sum_{k,j=1}^d \int_0^1 \E_U ( \partial_k \partial_j \cal{Q}_n(F_\l-y) H_{j,U}(G,1) (F-G)_k) d\l \\
&\quad = \sum_{k,j=1}^d \int_0^1 \E_U ( \partial_j \cal{Q}_n(F_\l-y) H_{k,U}( F_\l, H_{j,U}(G,1) (F-G)_k)) d\l  \\
&\quad = \sum_{k,j=1}^d \int_0^1 \E_U ( \partial_j \cal{Q}_n(F_\l-y) H_{k,U}(F_\l,  V_{j,k})) d\l 
\end{split}
\]
Now, applying first \eqref{wei1} and then \eqref{wei2}, with some computations in the same fashion as before, it is possible to show
\[
\begin{split}
& \| (H_{k,U}(F_\l , V_{j,k}))_{j=1,\dots,n} \|_{2n,U} \\
& \quad \leq C 
m_U(32n )^2 \Gamma_{F,G,U}(32 n)^5 n_{F,G,U}(32 n)^8 
\|F-G\|_{1,32n,U}.
\end{split}
\]
From \eqref{estQ} and H\"older as before,
\[
|J|\leq 
C
\left[
m_U(32n ) \Gamma_{F,G,U}(32 n)^2 n_{F,G,U}(32 n)^3
\right]^{4n^2}\|F-G\|_{1,32n,U}.
\]
The statement follows.
\epr

\bl{lemma2}
Let $U$ be a localizing r.v.,  $F=(F_1,\dots,F_n), G=(G_1,\dots,G_n)\in \DD^{3,\infty}$. If $\Gamma_{F,U}(q)<\infty$, $\Gamma_{G,U}(q)<\infty$ for any $q>1$,
there exists a constant $C$ depending only on the dimension $n$ such that
\[
\begin{split}
&|p_{F,U}(y)-p_{G,U}(y)|\\
&\, \leq C 
\left[ m_U(32n ) (\Gamma_{F,U}\vee\Gamma_{G,U})(32 n) n_{F,G,U}(32 n)\right]^{24 n^2}
\|\Delta_2(F,G)\|_{32 n,U}
\end{split}
\]
\el

\bpr
We denote in this proof $M=\hat{\g}_G (\g_{F_\l}-\g_G)$, and
define, as in \eqref{defU},  
\be{defV}
V=\prod_{1\leq i,j \leq n} \psi_{1/(8 n^2)}(M_{i,j}).
\ee
We have from Lemma \ref{lemma1} that if $\Gamma_{F,G,UV}(q)$ is finite for $q>0$ \be{den1}
\begin{split}
&|p_{F,UV}(y)-p_{G,UV}(y)|\\
&\,\leq C 
\left[ m_{UV}(32n ) \Gamma_{F,G,UV}(32 n) n_{F,G,UV}(32 n)\right]^{12 n^2}
\|\Delta_2(F,G)\|_{32 n,UV}
\end{split}
\ee
Remark
\[
\hat{\g}_G-\hat{\g}_{F_\l}=\hat{\g}_G (\g_{F_\l}-\g_G) \hat{\g}_{F_\l},
\]
so
\[
\|\hat{\g}_{F_\l}-\hat{\g}_G\|_{Fr} \leq
\| \hat{\g}_G (\g_{F_\l}-\g_G)\|_{Fr} \|\hat{\g}_{F_\l}\|_{Fr}
\]
On $V\neq 0$ we have $\| \hat{\g}_G (\g_{F_\l}-\g_G) \|_{Fr} \leq 1/2$, because of definition \eqref{defV}, so
\[
\|\hat{\g}_{F_\l}\|_{Fr} \leq 2 \|\hat{\g}_G\|_{Fr}
\]
and therefore 
\be{FGG}
\G_{F,G,UV}(32 n)\leq 2\G_{G,UV}(32 n)\leq 2\G_{G,U}(32 n).
\ee
Now, using \eqref{changewei} with $G=1$,
\[\begin{split}
p_{F,U (1-V)}(y)
&=\E_{U(1-V)}[\nabla Q(F-y),H_{U(1-V)}(F,1)]\\
&=\E_{U}[\nabla Q(F-y),(1-V)H_{U(1-V)}(F,1)]\\
&=\E_{U}[\nabla Q(F-y),H_{U}(F,1-V)] 
\end{split}\]
which implies, using as before \eqref{wei1} and \eqref{estQ}
\[
\begin{split}
&p_{F,U(1-V)}(y)
=E_{U(1-V)}\langle \nabla \cal{Q}_d(F-y),H_{U}(F,1-V)\rangle\\
&\quad
\leq C \left[
m_U(14n) \Gamma_{F,U}(14n)^2 n_{F,U}(14n)^3
\right]^{4n(n-1)} \|H_{U}(F,1-V)\|_{2n,U}\\
&\quad
\leq C \left[
m_U(24n) \Gamma_{F,U}(24n)^2 n_{F,U}(24n)^3
\right]^{8n(n-1)+1} \|1-V\|_{1,4n,U}
\end{split}
\]
and, using \eqref{inequality2},
\[
\|1-V\|_{1,4n,U} 
\leq C\|\hat{\g}_G (\g_{F_\l}-\g_G)\|_{1,4n,U}
\]
Now, we first apply H\"older inequality and then
\[
|\g_{F_\l}-\g_G| \leq C  |D ({F_\l}-G)| \,| D ({F_\l}+G)|
\]
and
\[
|D(\g_{F_\l}-\g_G) | \leq C  \left(|D({F_\l}-G)|+|D^{(2)} ({F_\l}-G)| \right) \left(|D({F_\l}+G)|+|D^{(2)} ({F_\l}+G)| \right)
\]
We find
\[
\|1-V\|_{1,4n,U} 
 \leq C\G_{G,U}(12n) n_{F,G,U}(12n) \|F-G\|_{2,12n,U}
\]
so
\[
p_{F,U(1-V)}(y)\leq C \left[
m_U(24n) (\Gamma_{F,U} \vee \Gamma_{G,U})(24n)^2 n_{F,G,U}(24n)^3
\right]^{8n^2} \|\Delta_2(F,G)\|_{32n,U}
\]
We conclude writing
\[
\begin{split}
|p_{F,U}(y)-p_{G,U}(y)|
&=
|p_{F,UV}(y)+p_{F,U(1-V)}(y)
-p_{G,UV}(y)-p_{G,U(1-V)}(y)|\\
&\leq
|p_{F,UV}(y)-p_{G,UV}(y)|+p_{F,U(1-V)}(y)+p_{G,U(1-V)}(y)
\end{split}
\]
and the statement follows easily.
\epr

\bpr \itshape (of Theorem \ref{genden}). \upshape  Let $V$ as in the last proof. We can write
\[
\begin{split}
p_{F,U}(y) & \geq p_{F,UV} (y)\geq p_{G,UV} (y) - |p_{F,UV} (y) - p_{G,UV} (y)|\\
&= p_{G,U} (y) - p_{G,U(1-V )}(y) - |p_{F,UV} (y) - p_{G,UV} (y)|.
\end{split}
\]
From \eqref{wei1} and \eqref{estQ} as before
\[
\begin{split}
p_{G,U(1-V)}(y)
& \leq C \left[
m_U(14n) \Gamma_{G,U}(14n)^2 n_{F,G,U}(14n)^3
\right]^{8n^2} \|\Delta_2(F,G)\|_{32n,U}.
\end{split}
\]
Using also \eqref{den1} and \eqref{FGG}
we obtain the desired lower bound for $p_F$.

For the upper bound we apply Proposition \ref{repden} localizing on $1-U$. We have
\[
p_{F,1-U}(x) =
\E_{(1-U)}\left[\nabla \Q_n(F-x) H_{(1-U)}(F,1)
\right]
=
\E\left[\nabla \Q_n(F-x) H_{(1-U)}(F,1) \,(1-U)
\right]\]
From \eqref{changewei}, $H(F,1-U)=(1-U) H_{(1-U)}(F,1)$, so
\[
p_{F,1-U}(x) =
\E [\nabla \Q_n(F-x) H(F,1-U)]
\]
Now we apply H\"older and find 
\[
p_{F,1-U}(x)=
\|\nabla\Q_n(F-x)\|_{\frac{2n}{2n-1}}   \| H(F,1-U)\|_{2n}
\]
We use \eqref{estQ}, with $U=1$, to deal with the gradient of the Poisson kernel:
\[
(\E|\nabla \cal{Q}_n(F-y)|^{2n/(2n-1)})^{(2n-1)/(2n)}
\leq C \left(
\Gamma_{F}(14n)^2 n_F(14n)^3
\right)^{4n(n-1)}.
\]
Now consider the non-localized version of \eqref{wei1}:
\[ 
\|H(F,V)\|_{p} \leq C \|V\|_{1,14n}
\Gamma_{F}(14n)^{2} n_{F}(14n)^{3} 
\]
and take $V=1-U$. We obtain
\be{denlocR}
p_{F,1-U}\leq C \|1-U\|_{1,14n}\left[
\Gamma_{F}(14n)^2 n_F(14n)^3
\right]^{4n^2}.
\ee
We apply now the lower bound result to $p_{G,U}$, interchanging the roles of $F$ and $G$, and find
\[
p_{F,U}(y)\leq p_{G,U}(y)
+ 
\left[ m_U(32n ) \Gamma_{F,U}(32 n) n_{F,G}(32 n)\right]^{24 n^2}
\|\Delta_2(F,G)\|_{32 n,U}.
\]
Putting together this inequality and \eqref{denlocR}, we have the upper bound. 
\epr

\subsection{Density estimates via local inversion}
\label{sectioninversefunction}

We recall some results from \cite{BCP1}. We see how to use the local inversion theorem to transfer a known estimate for a Gaussian random variable to its image via a function $\eta$ such that
\begin{equation*}
\eta\in C^3(\R^n,\R^n),\quad \eta(0)=0,\quad \l^*(\nabla\eta(0))\leq \frac{1}{2}.
\end{equation*}
Define, for $h>0$,
$$
c_*(\eta,h)=\sup_{|x|\leq 2h} \max_{i,j} |\partial_i\eta^j(x)|
$$
and
\[
c_2(\eta) =\max_{i,j=1,..,n}\sup_{|x|\leq 1}|\partial^2_{ij}\eta(x)|,\quad
c_3(\eta) =\max_{i,j,k=1,..,n}\sup_{|x|\leq 1}|\partial^3_{ijk}\eta(x)|,
\]
Let now $\Theta$ be a $n$-dimensional centered Gaussian variable with covariance matrix $Q$, non-degenerate. Denote by $\underline{\l}$ and $\overline{\l}$ the lower and the upper eigenvalues of $Q$. 
Suppose to have $r>0$ such that
\be{hpimpl}
c_*(\eta,16 r)\leq \frac{1}{2n} \sqrt{\frac{\underline{\l}}{\overline{\l}}},\quad \quad
r\leq h_\eta =\frac{1}{16 n^2 (c_2(\eta)+\sqrt{c_3(\eta)})}.
\ee
We take a localizing function as in \eqref{defU}: $U_r=\prod_{i=1}^n \psi_{r}(\Theta_i)$. We also define $\Phi(\theta)=\theta+\eta(\theta)$. 
\bl{invfun2}
The density $p_{G,U_r}$ of
\[
G:=\Phi(\Theta)=\Theta+\eta(\Theta)
\]
under $\PR_{U_r}$ has the following bounds on $B(0,r)$:
\[
\frac{1}{C \det Q^{1/2}
}\exp\left(-\frac{C}{\underline{\l}}|z|^2 \right) \leq 
p_{G,U_r}(z) \leq \frac{C}{\det Q^{1/2} }\exp\left(-\frac{1}{C \overline{\l}}|z|^2 \right)
\]
\el
This result is proved in \cite{BCP1} under a slightly stronger constraint on $r$, but going trough the proof it is easy to see that what we suppose here is enough. For details see \cite{tesipigato}. The proof is quite standard and follows from the local inversion theorem (see \cite{rudin-principles} for a standard version of this theorem).

\section{Density estimates of the diffusion process}\label{sshortime}
In this section we prove lower and upper bounds for the density of $X_\delta$. 

\subsection{Development}
In this section, in order to lighten the notation, we do not mention the dependence of the parameters on the initial condition (so, for example, we write $A$ instead of $A(x_0)$, and so on). We need to introduce some notation. Consider a small time $\delta\in(0,1]$. We define
\bi
\item The translated initial condition
\[
\hat{x}_0=x_0+b(x_0)\delta.
\]
\item The matrices $\bar{A}$ and $\bar{A}_\delta$ as
\[
\bar{A}=(\s+\d\partial_b\s,[\s,b])
\]
and
\[
\bar{A}_\delta=\left(\delta^{1/2} (\s+\d\partial_b\s), \delta^{3/2} [\s,b]\right).
\]
Recall \eqref{defA}, \eqref{defAR}, and remark that {\bf A1} implies that these matrices are always invertible if $\d$ is small enough.
\item The Gaussian r.v.
\[
\Theta=\left( \begin{array}{l}
\Theta_1 \\
\Theta_2
\end{array}\right)
=\left( \begin{array}{l}
\delta^{-1/2} W_\delta\\
\delta^{-3/2} \int_0^\delta (\delta-s)dW_s
\end{array}\right).
\]
\item The polynomial of degree 3 and direction $\s(x_0)$ (recall $\kappa_\s$ defined in \eqref{lambdasigma}):
\be{defeta}
\eta(u)= \left( \frac{\kappa_\s(x_0) }{2} u^2 + 
\frac{(\partial_\s \kappa_\s+\kappa_\s^2)(x_0)}{6} u^3 \right)\s(x_0).
\ee
\item The principal term
\be{defG}
G=\Theta+\tilde{\eta}_\d(\Theta)
\ee
where $\tilde{\eta}_\d(\Theta) = \bar{A}^{-1}_\delta \eta(\delta^{1/2}\Theta_1)$.

\item The remainder $R_\d$:
\be{remainder}
\begin{split}
R_\d=&\int_0^\d\int_0^s \left(\partial_b\sigma(X_u)-\partial_b\sigma(x_0)\right) du \circ dW_s\\
&+\int_0^\d\int_0^s \left(\partial_\sigma b(X_u)-\partial_\sigma b(x_0) \right) \circ dW_u ds\\
&+\int_0^\d \int_0^s \partial_b b(X_u) du ds\\
&+\int_0^\d\int_0^s \int_0^u \left(\partial_\sigma \partial_\sigma \sigma(X_v)-\partial_\sigma \partial_\sigma \sigma(x_0)\right) \circ dW_v\circ dW_u\circ dW_s\\
&+\int_0^\d\int_0^s\int_0^u \partial_b \partial_\sigma \sigma(X_v) \circ dv \circ dW_u\circ dW_s.
\end{split}
\ee
Notice that $R_\d \sim \cal{O}(\d^2)$. 
We also denote $\tilde{R}_\delta:=\bar{A}_\delta^{-1} R_\delta$.
\ei
We now prove that the following decomposition holds:
\be{dec}
X_\delta=\hat{x}_0+\bar{A}_\delta (G+\tilde{R}_\delta)
\ee
This is a main tool in our approach. Indeed, we find Gaussian bounds for the density of the variable $F:=\bar{A}_\d^{-1}(X_\d-\hat{x}_0)=G+\tilde{R}_\d$ in the Euclidean metric of $\R^2$. The fact that in Theorem \ref{mtstime} the bounds for the diffusion are in the $A_\d(x_0)$-norm follows from the change of variable suggested by $\eqref{dec}$.
 
Let us prove \eqref{dec}. We write the stochastic Taylor development of $X_t$ with a remainder of order $t^2$:
$$
X_t=x_0+b(x_0)t+U_t+R_t,
$$
where
\begin{align*}
U_t&= \sigma(x_0)W_t+\partial_{\sigma} \sigma(x_0) \int_0^t W_s\circ dW_s \\
&+ \partial_{\sigma}\partial_{\sigma} \sigma(x_0) 
\int_0^t\int_0^s W_u\circ dW_u\circ dW_s\\
&+\partial_b \sigma(x_0)\int_0^t  s dW_s
+\partial_{\sigma} b(x_0)\int_0^t W_s ds
\end{align*}
Now we write
\begin{align*}
&\int_0^t W_s ds = \int_0^t (t-s)dW_s\\
&\int_0^t  s dW_s= - \int_0^t (t-s) dW_s + t W_t
\end{align*}
Therefore
\begin{align*}
U_t&= (\sigma(x_0) + t\partial_b \sigma(x_0))W_t
+(\partial_{\sigma} b(x_0)-\partial_b \sigma(x_0))\int_0^t (t-s)dW_s\\
&+\partial_{\sigma} \sigma(x_0) \frac{W_t^2}{2}+\partial_{\sigma}\partial_{\sigma} \sigma(x_0) 
\frac{W_t^3}{6}
\end{align*}
So we have the following decomposition of $X_t$:
\be{devn}
X_t=x_0+b(x_0)t+(\sigma(x_0)+t\partial_b\sigma(x_0))W_t+[\s,b](x_0)\int_0^t(t-s)dW_s+\eta(W_t)+R_t
\ee
where $x_0$ is the initial condition. 
Remark that {\bf A3} implies that both the coefficients of $\eta$ have the same direction as $\s(x_0)$:
\[
\eta(u)=\frac{\partial_\sigma \sigma(x_0)}{2} u^2 + 
\frac{\partial_\sigma \partial_\sigma\sigma(x_0)}{6} u^3=
\left( \frac{\kappa_\s (x_0)}{2} u^2 + 
\frac{(\partial_\s \kappa_\s+\kappa_\s^2)(x_0)}{6} u^3 \right)\s(x_0).
\]

\subsection{Preliminary estimates}\label{prel}

We introduce the following class of constants:
\be{def:constant}
\cal{C}=\left\{C>0: C= K\left(\frac{\rho}{\l_*(A(x_0))}\right)^q,\exists K,q\geq 1 \right\}
\ee
We stress that the constants defined above depend on the parameters of the diffusion through the ratio $\rho/\l_*(A(x_0))$ (cf. \textbf{A1}, \textbf{A2}), but $K,q$ do not depend on $\s,b$. 
We will also denote by $1/\cal{C}=\{\d>0\,:1/\d\in \cal{C}\}$. 

We keep using the notations of the previous development.
\bl{lemma41}
There exist $L_1,L_2,K_1,K_2$ positive constants not depending on the parameters, $\d^*\in 1/\cal{C}$ such that: for any fixed $r>0$ and $\d$ such that $\d\leq \d^* \exp\left(-2 L_1 r^2\right)$, let $G=\Theta+\tilde{\eta}_\d(\Theta)$ be the r.v. defined in \eqref{defG}; let $U_r$ be the localizing r.v. defined in \eqref{defU}, and $p_{G,U_r}$ the local density of $G$; then the following estimate holds for $|z|\leq r$:
\be{ineqstep1}
K_1 \exp\left(-L_1 |z|^2\right) 
\leq p_{G,U_r}(z)\leq 
K_2 \exp\left(-L_2 |z|^2\right).
\ee
\el

\bpr 
In what follows, $C\in \cal{C}$, and may vary from line to line (meaning that $K,q$ may vary in \eqref{def:constant}). We start by computing the derivatives of $\eta$:
\begin{align*}
\eta(y)&=\left( \frac{\kappa_\s}{2} y^2 + 
\frac{\partial_\s \kappa_\s +\kappa_\s^2}{6} y^3 \right) \s \\
\eta '(y)&=\left( \kappa_\s y + 
\frac{\partial_\s \kappa_\s +\kappa_\s^2}{2} y^2 \right) \s \\
\eta ''(y)&=( \kappa_\s  + (\partial_\s \kappa_\s +\kappa_\s^2) y)\s\\
\eta '''(y)&= (\partial_\s \kappa_\s +\kappa_\s^2)\s.
\end{align*}
By the definition of $\bar{A}_\delta^{-1}$, 
\begin{align*}
\bar{A}_\delta^{-1} \delta^{1/2}(\sigma+\delta\partial_b \sigma)=(1,\,0)^T.
\end{align*}
Therefore
\begin{align*}
\bar{A}_\delta^{-1} \sigma=  
\delta^{-1/2}(1,\,0)^T-\bar{A}_\delta^{-1} \delta\partial_b \sigma.
\end{align*}
By  \eqref{basicmatrnorm2} and \eqref{eq3} (see the appendix) we have $|\bar{A}_\delta^{-1} \delta\partial_b \sigma|\leq C \delta^{-1/2}$, so that $|\bar{A}_\delta^{-1} \sigma|
\leq C \delta^{-1/2}$. We stress that this upper bound is $\d^{-1/2}$ in contrast with $\d^{-3/2}$ in \eqref{basicmatrnorm2}, because $\bar{A}_\d$ works in the specific direction $\s$.
Now we can estimate the norms of $\tilde{\eta}_\d$ and its derivatives. Since they are collinear with $\s$, we have
\[
\begin{split}
|\tilde{\eta}_\d(u)|&= |\bar{A}_\delta^{-1} \, \eta(\delta^{1/2} u_1)|
\leq C (|u_1|^2 \delta^{1/2}+ |u_1|^3\d )\\
|\partial_{u_1}\tilde{\eta}_\d(u)|&= |\bar{A}_\delta^{-1} \delta^{1/2}\eta '(\delta^{1/2} u_1)| \leq C (|u_1|\d^{1/2} + |u_1|^2 \d)\\
|\partial^2_{u_1}\tilde{\eta}_\d (u)|&= |\bar{A}_\delta^{-1} \delta\eta ''(\delta^{1/2} u_1)| \leq C (\delta^{1/2}  +|u_1| \d)\\
|\partial^3_{u_1}\tilde{\eta}_\d(u)|&= |\bar{A}_\delta^{-1} \delta^{3/2}\eta '''(\delta^{1/2} u_1)|  \leq C\d\\
|\partial_{u_2}\tilde{\eta}_\d(u)|&=0.
\end{split}
\]
So, referring to the notation of Section \ref{sectioninversefunction}, we have
\be{dereta}
\begin{split}
c_*(\tilde{\eta}_\d,h)&=\sup_{|u|\leq 2h} \max_{i,j} \left|\partial_i\tilde{\eta}_\d^j(u) \right| \leq C h \d^{1/2}\\
c_2(\tilde{\eta}_\d)&=\max_{i,j} \sup_{|u|\leq 1}  \left|\partial^2_{i,j}\tilde{\eta}_\d(u) \right|\leq C \d^{1/2} \\
c_3(\tilde{\eta}_\d)&=\max_{i,j,k} \sup_{|u|\leq 1}  \left|\partial^3_{i,j,k}\tilde{\eta}_\d(u) \right|\leq C \d.
\end{split}
\ee
We first want to apply Lemma \ref{invfun2} to $G=\Theta+\tilde{\eta}_\d(\Theta)$.
Here $n=2$, and the covariance matrix of $\Th$ is
\[
\g_\Theta=\left(\begin{array}{cc}
1 & 1/2 \\
1/2 & 1/3
\end{array}\right).
\] 
It has 2 positive eigenvalues, $0<\l_1<\l_2$, and $\det(\g_\Th)=1/12$. 
We are supposing here $ \d\leq \d^* \exp\left(-2 L_1 r^2\right)\leq 
\d^* /r^2$. Since
\[
h_{\tilde{\eta}_\d}=\frac{1}{64 (c_2(\tilde{\eta}_\d)+\sqrt{c_3(\tilde{\eta}_\d)})}\geq \frac{1}{C_1 \sqrt{\d}} \geq  \frac{r}{C_1 \sqrt{\d^*}}
\]
and
\[
c_*(\tilde{\eta}_\d,16r)\leq C_2 r \sqrt{\d} \leq C_2  \sqrt{\d^*},
\]
choosing $\d^*\leq \frac{1}{16} \frac{\l_1}{\l_2} \frac{1}{C_1^2 C_2^2}$ the conditions \eqref{hpimpl} are satisfied: 
\be{ara}
c_*(\tilde{\eta}_\d,16 r)\leq \frac{1}{4} \sqrt{\frac{\l_1}{\l_2}},\quad \quad
r\leq h_{\tilde{\eta}_\d}
\ee
So there exist $L_1,L_2,K_1,K_2$ universal constants, such that for $|z|\leq r$,
\[
K_1 \exp\left(-L_1 |z|^2\right) 
\leq p_{G,U_r}(z)\leq 
K_2 \exp\left(-L_2 |z|^2\right).
\]
\epr

The following lemma is a slight modification of Lemma 2.3.1. in \cite{Nualart:06}.
\bl{lemma:matrixmoment}
Let $\g$ be a symmetric non-negative definite $n\times n$ matrix. We assume that, for fixed $p\geq 2$, $\E\left[ \|\g\|_{Fr}^{p+1} \right]< \infty$, and that $\exists\,\ve_0>0$ s.t. for $\ve \leq \ve_0$,
\[
\sup_{|\xi|=1} \PR [\langle\g \xi,\xi \rangle<\ve]\leq \ve^{p+2n}
\]
Then there exist a constant $C$ depending only on the dimension $n$ such that
\[
\E \left[ \l_*(\g)^{-p}\right] \leq C \E \left[ \|\g\|_{Fr}^{p+1}\right] \ve_0^{-p}.
\]
\el
We consider now 
\be{defF}
F=\bar{A}_\d^{-1}(X_\d-\hat{x}_0).
\ee
We will use the general estimates of section \ref{sec:malliavin}. We denote by $D$ the Malliavin derivative with respect to $W$, the Brownian motion driving \eqref{eqn}.
We first prove that the moments of $\l^*(\g^{-1}_F)=\l_*(\g_F)^{-1}$ are bounded, and these bounds do not depend on $\d$. This result looks interesting by itself, since it means that we are able to account precisely of the scaling of the diffusion in the two main directions $\s$ and $[\s,b]$. In this particular case this is a refinement of the classical result on the bounds of the Malliavin covariance under the (weak) H\"{o}rmander condition (cf. \cite{Nualart:06}, \cite{KusuokaStroock:85}, \cite{Norris:86}). 

\bl{n}
Let $F=\bar{A}_\d^{-1}(X_\d-\hat{x}_0)$. For any $p>1$, there exists $C\in \cal{C}$ such that for any $\d\leq1$,
$\Gamma_F(p)\leq e^C$.
\el
\bpr
Following \cite{Nualart:06} we define the tangent flow of $X$ as the derivative with respect to the initial condition of $X$, $Y_t:=\partial_x X_t$. We also denote its inverse $Z_t=Y_t^{-1}$. They satisfy the following stochastic differential equations
\[
\begin{split}
Y_t&=Id + \int_0^t \nabla\s(X_s) Y_s \circ dW_s+\int_0^t \nabla b(X_s) Y_s ds\\
Z_t&=Id - \int_0^t Z_s \nabla\s(X_s) \circ dW_s-\int_0^t Z_s \nabla b(X_s)ds
\end{split}
\]
The Malliavin derivative of $X$ is
\[
D_s X_t=Y_t Z_s \s(X_s),
\]
so
\[
D_s F = D_s \bar{A}_\d^{-1} (X_\d-\hat{x}_0) = \bar{A}_\d^{-1} Y_\d Z_s \sigma(X_s).
\]
We define
\[
\bar{\g}_\d=\int_0^\d A_\d^{-1} Z_s \sigma(X_s) \sigma(X_s)^T Z_s^T A_\d^{-1,T} ds.
\]
Then
\[
\g_F=\langle D F,D F\rangle= \bar{A}_\d^{-1} Y_\d A_\d \bar{\g}_\d A_\d^T Y_\d^T \bar{A}_\d^{-1,T}.
\]
Remark that
\[
\g_F^{-1}= 
\bar{A}_\d^T Z_\d^T  A_\d^{-1,T}
\bar{\g}_\d^{-1} 
A_\d^{-1} Z_\d \bar{A}_\d,
\]
and that in this representation we have both $A_\d$ and its ``perturbed" version $\bar{A}_\d$. We have to check the integrability of 
$\l_*(\g_F)^{-1}=\l^*(\g_F^{-1})$. Recall that $\l^*(\cdot)$ is a norm on the set of matrices, and that for two $2\times 2$ matrices $M_1,M_2$,  $\l^*(M_1 M_2)\leq 2 \l^*(M_1) \l^*(M_2) $. We have
\[
\l_*(\g_F)^{-1} \leq 4\l^*( \bar{\g}_\d^{-1}) \l^*(A_\d^{-1} Z_\d \bar{A}_\d)^2,
\]
We need to bound
$A_\d^{-1} Z_\d \bar{A}_\d$, which we expect to be close to the identity matrix for small $\d$, and $\bar{\g}_\d^{-1} $.

We take care first of the moments of $\l^*( \bar{\g}_\d^{-1})$.
We use the following representation, holding for general $\phi$, which follows applying Ito's formula (details in \cite{Nualart:06})
\be{devZphi}
Z_t \phi(X_t)=
\phi(x_0)+\int_0^t Z_s [\s,\phi](X_s) dW_s^k+
\int_0^t Z_s \left\{[b,\phi]+\frac{1}{2}[\s,[\s,\phi]] \right\}(X_s)ds
\ee
Taking $\phi=\s$ the representation above reduces to 
\be{zl}
\begin{split}
Z_t \s(X_t) &=\s(x_0)+\int_0^t Z_s [b,\s](X_s)ds\\
&=\s(x_0)+  t [b,\s](x_0) + L_t,
\end{split}
\ee
with $L_t=\int_0^t Z_s [b,\s](X_s) - Z_0 [b,\s](x_0)ds$. Notice that $L_t\sim \cal{O}(t^{3/2})$. Using \textbf{A2} one gets
\[
\E\left[
\l^*\left( \int_0^{\d \ve} L_s L_s^T ds \right)^q
\right] \leq
\E\left[
\left\|\int_0^{\d \ve} L_s L_s^T ds \right\|_{Fr}^q
\right] \leq  e^{C'} (\d \ve)^{4q}, \quad \forall q>0, \quad \exists C'\in \cal{C}
\]
($e^{C'}$ comes from Gronwall inequality).
We have
\[
\begin{split}
A_\d^{-1} Z_s \s(X_s)& =
A_\d^{-1}(\s(x_0)+  s [b,\s](x_0) + L_s)\\
&= \frac{1}{\d^{1/2}} 
\left( \begin{array}{c}
    1\\ -s/\d
\end{array} \right)+
A_\d^{-1} L_s
\end{split}
\]
For constant $c$ and fixed $\ve$, we introduce the stopping time
\[
S_\ve = \inf  \left\{ s \geq 0: \l^*\left( \int_0^s L_u L_u^T du \right) \geq c(\d\ve)^3  \right\} \wedge \d,
\]
We have
\[
\l^*\left( A_\d^{-1} \int_0^{S_\ve} L_u L_u^T du A_\d^{-1,T} \right)
\leq 4
\l^*\left( A_\d^{-1}\right)^2 \l^*\left(\int_0^{S_\ve} L_u L_u^T du  \right)
\leq \frac{C''}{\d^3}c(\d\ve)^3 
\]
where $C''\in \cal{C}$.  We fix $c=\frac{1}{64 C''}$, so
\be{boundstoppingtime}
\l^*\left( A_\d^{-1} \int_0^{S_\ve} L_u L_u^T du A_\d^{-1,T} \right)
\leq \frac{\ve^3}{64}
\ee
Now we suppose to be on the event $\{\frac{S_\ve}{\d}\geq \ve\}$. 
Applying first inequality 
\[
\langle (v+R)(v+R)^T\xi,\xi\rangle \geq \frac{1}{2}\langle vv^T\xi,\xi\rangle-\langle RR^T\xi,\xi \rangle,
\]
which holds for any vectors $v,\, R,\,\xi$, and then \eqref{zl}, we obtain
\[
\begin{split}
\bar{\g}_\d &=\int_0^\d A_\d^{-1} Z_s \sigma(X_s) \sigma(X_s)^T Z_s^T A_\d^{-1,T} ds\\
&\geq \int_0^{S_\ve} A_\d^{-1} Z_s \sigma(X_s) \sigma(X_s)^T Z_s^T A_\d^{-1,T} ds\\
&\geq\frac{1}{2} \int_0^{S_\ve} \frac{1}{\d} 
\left( \begin{array}{cc}
    1 & -s/\d \\ -s/\d & (s/\d)^2
\end{array} \right) ds
- A_\d^{-1} \int_0^{S_\ve} L_s L_s^T ds A_\d^{-1,T}.
\end{split}
\]
We have
\[
\int_0^{S_\ve} \frac{1}{\d} 
\left( \begin{array}{cc}
    1 & -s/\d \\ -s/\d & (s/\d)^2
\end{array} \right) ds
\geq
\int_0^{\d \ve} \frac{1}{\d} 
\left( \begin{array}{cc}
    1 & -s/\d \\ -s/\d & (s/\d)^2
\end{array} \right) ds
\geq
\left( \begin{array}{cc}
\ve & -\frac{\ve^2}{2} \\ -\frac{\ve^2}{2} & \frac{\ve^3}{3}
\end{array} \right)
\geq
Id_2\,\frac{\ve^3}{16},
\]
so, from \eqref{boundstoppingtime},
\[
\langle\bar{\g}_\d \xi, \xi\rangle \geq \frac{1}{2}\frac{\ve^3}{16} |\xi|^2 - \frac{\ve^3}{64} |\xi|^2=\frac{\ve^3}{64} |\xi|^2, \quad \forall |\xi|=1.
\]
Now, remark that $t\rightarrow \l^*\left(\int_0^t L_s L_s^T ds \right)$ is increasing. For any $q>0$
\[
\begin{split}
\PR(S_\ve < \d\ve)
&\leq
\PR\left(
\l^*\left(\int_0^{\d \ve} L_s L_s^T ds \right)^q
\geq c^q(\d\ve)^{3q}
\right)\\
&\leq
\frac{\E\left[
\l^*\left(\int_0^{\d \ve} L_s L_s^T ds \right)^q
\right]}
{c^q(\d\ve)^{3 q}}\\
&\leq
\frac{e^{C'} (\d \ve)^{4q}}
{c^q(\d\ve)^{3 q}}\leq
\frac{ e^{C'}}{c^q} (\d \ve)^q \leq \ve^{q/2} 
\end{split}
\]
for $\d\leq 1$, for $\ve\leq \ve_0=e^{-C'''}$ with $C'''\in\cal{C}$.
Therefore, for any $q$, for any $\ve\leq\ve_0$, $\d\leq 1$,
\[
\PR(\langle \bar{\g}_\d\xi, \xi\rangle < \ve^{3}/64) \leq
\PR[S_\ve < \d\ve]
 \leq \ve^{q/2}
\]
Now we apply Lemma \ref{lemma:matrixmoment}. We obtain
\[
\E \l^*( \bar{\g}_\d^{-1})^q=\E \l_*( \bar{\g}_\d)^{-q}\leq e^C
\]
for $\d\leq 1$, $C\in \cal{C}$.

We consider now $A_\d^{-1} Z_\d \bar{A}_\d$. Applying \eqref{devZphi} and {\bf A3}, one can prove that
\[
Z_t \s(x_0)= (1-\kappa_\s(x_0)W_t)\s(x_0) + J_t,
\]
with $J_t\sim\cal{O} ( t)$. So
\[
Z_\d \bar{A}_\d = \left( 
\sqrt{\d}(1-\kappa_\s(x_0)W_\d)\s(x_0),\, 0 
\right)+ M_\d
\] 
where $M_\d$ is a $2\times 2$ matrix with $\E\l^*(M_\d)^q\leq e^C\d^{3q/2}$, $C\in \cal{C}$. This estimate follows again from \textbf{A2}.
 Since $A_\d=(\d^{1/2}\s(x_0), \d^{3/2}[\s,b](x_0))$
\[
A_\d^{-1} \left( 
\sqrt{\d}(1-\kappa_\s(x_0)W_\d)\s(x_0),\, 0 
\right)
=
\left( \begin{array}{cc}
1-\kappa_\s(x_0)W_\d  & 0  \\ 0  & 0
\end{array} \right)
\] 
and $\E |1-\kappa_\s(x_0)W_\d|^q\leq C\in \cal{C}$. 
Clearly $\E \l^*\left(A_\d^{-1} M_\d \right)^q\leq e^C$, $C\in \cal{C}$, so 
\[
\E\l^*(A_\d^{-1} Z_\d \bar{A}_\d)^q\leq e^C, \quad C\in \cal{C}.
\]
\epr

\subsection{Two-sided bound for the density of $X_\d$}

In this section we prove the short time density estimate \eqref{stintro}. We start with the following lemma, which is a density estimate for the ``renormalized'' random variable $F$ (see \eqref{defF}). 
We use Theorem \ref{genden} to recover estimates for $p_F$ from \eqref{ineqstep1}. We will need the preliminary estimates of Section \ref{prel}.

\bl{rescdiff}
Recall \eqref{def:constant}, the definition of $\C$, and that, for fixed $\d>0$, we set $F=\bar{A}_\d^{-1}(X_\d-\hat{x}_0)$ and $p_F$ is its density.
\begin{enumerate}
\item
There exist $C,C^*,L\in \cal{C}$ such that the following holds. We set $\d^*= e^{-C^*}$. For any fixed $r>0$, if $\d\leq \d^* \exp\left(-L r^2\right)$, for $|z|\leq r$ we have
\[
\frac{1}{C}\exp\left(-C |z|^2\right) 
\leq p_{F}(z)
\]
\item There exists $\d^*\in 1/\C;\,C,L\in \cal{C}$ such that: 
for any fixed $r>0$, if $\d\leq \d^*\exp\left(-L r^2\right)$, 
for $|z|\leq r$, we have
\[
p_{F}(z)\leq 
e^C \exp\left(-C^{-1} |z|^2\right).
\]
\end{enumerate}
\el

\bpr 
We apply Theorem \ref{genden}. Here $n=2$, so $32n=64$.\\
\ben 
\item (lower bound) 
Let $L_1$ be the constant in Lemma \ref{lemma41}. We first prove the lower bound for $r\geq\frac{1}{\sqrt{L_1}}=:\tilde{r}$.

We start checking that $C_1$ in Remark \ref{c} is in $\cal{C}$.
From\eqref{normlocalization} and $r\geq\frac{1}{\sqrt{L_1}}$,
\[
m_{U_r}(64)
\leq  
C\left( 1+\frac{\|\Theta\|^2_{2,64}}{r^2}\right)\leq  
C\in \cal{C}.
\]
Recall that $G=\Th+\tilde{\eta}_\d(\Th)$, where $\Th$ is a Gaussian with covariance (and also Malliavin covariance matrix) given by
\[
\g_\Theta=\left(\begin{array}{cc}
1 & 1/2 \\
1/2 & 1/3
\end{array}\right).
\] 
This matrix has 2 positive eigenvalues, $0<\l_1<\l_2$. 
Recall also that the Malliavin derivative $D$ is taken with respect to the Brownian motion $W$ driving \eqref{eqn}.
We consider $\Gamma_{G,U_r}=1+(\E_{U_r} \l_*(\g_G)^{-p})^{1/p}$. 
\[
\begin{split}
\langle\g_G \xi, \xi \rangle & = \int_0^\d \langle D_s G,\xi \rangle^2\\
&\geq \int_0^\d \frac{1}{2}\langle D_s \Th,\xi \rangle^2 - \langle D_s \tilde{\eta}_\d (\Th),\xi \rangle^2 ds\\
&=S_1+S_2.
\end{split}
\]
We have
\[
S_2= \int_0^\d \langle \nabla\tilde{\eta}_\d(\Th) D_s \Th,\xi \rangle^2 ds
=  \int_0^\d \langle D_s \Th, \nabla\tilde{\eta}_\d(\Th)^T \xi \rangle^2 ds
\leq \l_2 \|\nabla\tilde{\eta}_\d(\Th)\|_{Fr}^2 |\xi|^2
\]
and $S_1\geq \l_1/2$, so
\[
\begin{split}
\l_*(\g_G) \geq \l_1 \left(\frac{1}{2} - \frac{\l_2}{\l_1} \|\nabla\tilde{\eta}_\d(\Th)\|_{Fr}^2\right).
\end{split}
\]
Recall $c_*(\tilde{\eta}_\d,h)=\sup_{|x|\leq 2h} \max_{i,j} |\partial_i\tilde{\eta}_\d^j(x)|$, so on the event $\{U_r\neq 0\}$ we have $|\Th|\leq 4r$ and 
$\|\nabla\tilde{\eta}_\d(\Th)\|_{Fr} \leq 2 c_*(\tilde{\eta}_\d, 16r)$. We proved in \eqref{ara} that $c_*(\tilde{\eta}_\d,16 r)\leq \frac{1}{4} \sqrt{\frac{\l_1}{\l_2}}$, so it  follows
\[
\|\nabla\tilde{\eta}_\d(\Th)\|_{Fr} \leq \frac{1}{2}\sqrt{\frac{\l_1}{\l_2}},
\]
and therefore $ \l_*(\g_G)\geq \l_1/4$, which implies $\Gamma_{G,U_r}(64)\leq C$.  
Recall \eqref{dec} and \eqref{defF}. Standard computations usign {\bf A2} and Gronwall lemma  give $n_{F,G,U_r}(p)\leq e^C,\,C\in \cal{C}$,
so from Theorem \ref{genden} we have that $\exists C\in \cal{C}$ such that for $|z|\leq r$
\[
p_F(z) 
\geq p_{G,U_r}(z)
- e^C \|\tilde{R}_\d\|_{64,U_r}\geq 
K_1 \exp\left(-L_1 |z|^2\right) 
- e^C \|\tilde{R}_\d\|_{64,U_r}.
\]
Recall \eqref{remainder}. By using {\bf A2}, one can show that $\| R_\d \|_{2,p}\leq e^{C} \d^2$, with $C\in \C$. 
So, from \eqref{basicmatrnorm1} with $\bar{A}_\d$ instead of $A_\d$,
\[
\|\tilde{R}_\d\|_{64,U_r}=\|\bar{A}_\d^{-1} R_\d\|_{64,U_r} \leq  e^C \d^2/\d^{3/2}  = e^C\sqrt{\d},
\]
so there exists $\bar{C}\in \cal{C}$ such that $p_F(z) \geq K_1 \exp\left(-L_1 |z|^2\right) - e^{\bar{C}}\sqrt{\d}$. We have that, for $r\geq \tilde{r}$,  if
\be{deltabound}
\d \leq \left(\frac{K_1 \exp(-\bar{C}) \exp(- L_1 r^2)}{2} \right)^2
=\left(\frac{K_1 \exp(-\bar{C}) }{2} \right)^2\, \exp(- 2L_1 r^2)
\ee
the following lower bound holds for $|z|\leq r$:
\[
p_F(z) \geq \frac{K_1}{2}\exp\left(-L_1 |z|^2\right)
\]
and this implies Lemma \ref{rescdiff}-(1) for $r\geq \tilde{r}$. We take now $0< r\leq \tilde{r}$. Remark that $\exp(- 2)=\exp(- 2L_1 \tilde{r}^2)$. We can suppose $\d^*\leq \left(\frac{K_1 \exp(-\bar{C}-1)}{2} \right)^2$, so
\[
\d \leq  \left(\frac{K_1 \exp(-\bar{C}-1) }{2} \right)^2=\left(\frac{K_1 \exp(-\bar{C}) }{2} \right)^2\exp(- 2L_1 \tilde{r}^2).
\]
If $|z|\leq r$, then $|z|\leq \tilde{r}$, and we apply what we have just proved for $r\geq \tilde{r}$, taking $\tilde{r}$ as radius. The following holds:
\[
p_F(z) \geq \frac{K_1}{2}\exp\left(-L_1 |z|^2\right).
\]

\item (upper bound). The proof of the upper bound follows again from Theorem \ref{genden}. We deal with $C_2$ exactly as for the lower bound, with the difference that we need a bound for $\Gamma_F(64)$ instead of $\Gamma_{G,U_r}(64)$. This is proved in  Lemma \ref{n}. 
As before, we first suppose $r\geq \frac{1}{\sqrt{L_2}}$, where $L_2$ is the constant in Lemma \ref{lemma41}.
We obtain
\[
p_F(z) 
\leq  
K_2 \exp\left(-L_2 |z|^2\right) 
+ e^{\bar{C}} (\sqrt{\d}+\|1-U_r\|_{1,28})
\]
$\bar{C}\in \cal{C}$. We fix $L\in \cal{C}$ and take $\d \leq \exp(- L r^2)$, 
and we also need to prove that $\|1-U_r\|_{1,28}$ decays as $C\exp(-{C^{-1}} |z|^2)$ for $|z|\leq r$. This follows from \eqref{inequality2}: $\exists C\in \cal{C}$ such that
\[
\|1-U_r\|_{1,28}\leq \sum_{i=1,2} \PR(|\Th_i|> r)^{\frac{1}{56}} C(1+1/r) \leq C e^{-C^{-1} r^2}. 
\]
We have the desired result for $r\geq \frac{1}{\sqrt{L_2}}$. Now, we take $r\leq \frac{1}{\sqrt{L_2}}$. If $|z|\leq r$, then $|z|\leq \frac{1}{\sqrt{L_2}}$, and we can apply the result already proved for $r\geq \frac{1}{\sqrt{L_2}}$, taking $\frac{1}{\sqrt{L_2}}$ as radius.. Then, we prove as in $(1)$ that the result can be extended to all $r>0$.
\een
\epr

We are now ready to prove the main theorem in short time.
\bt{mtstime}
Suppose {\bf A1}, {\bf A2}, {\bf A3} hold. Let $(X_t)_{t\in [0,T]}$ be the solution of $\eqref{eqn}$, and for $t\in [0,T]$, let $p_t(x_0,y)$ be the density of $X_t$ at $y$. 
\begin{enumerate}
\item
There exist $C,C^*,L\in \cal{C}$ such that the following holds. We set $\d^*= e^{-C^*}$. For any fixed $r>0$, if $0<\d\leq \d^* \exp\left(-L r^2\right)$, setting $\hat{x}_0=x_0+b(x_0)\d$, for $|y-\hat{x}_0|_{A_\d(x_0)} \leq r$ we have
\[
\frac{1}{C\d^2}\exp\left(-C |y-\hat{x}_0|_{A_\d(x_0)} ^2\right) 
\leq p_\d(x_0,y)
\]
\item There exists
$\d^*\in 1/\C,\,L,C\in \cal{C}$ such that: 
for any fixed $r>0$, if $0<\d\leq \d^*\exp\left(-L r^2\right)$, setting $\hat{x}_0=x_0+b(x_0)\d$, 
for $|y-\hat{x}_0|_{A_\d(x_0)}\leq r$, we have
\[
p_\d(x_0,y)\leq 
\frac{e^C}{\d^2}
\exp\left(-C^{-1} |y-\hat{x}_0|_{A_\d(x_0)}^2 \right).
\]
\end{enumerate}
\et
\bpr
We write the expectation of $f(X_\d)$ for a function $f$ with support included in $B(0,r)$. We  get
\[
\E[f(X_\d)]  = \E[f(\hat{x}_0+ \bar{A}_\d F ) ]
= \int f(\hat{x}_0+ \bar{A}_\d z) p_F(z) dz.
\]
With $\d,r$ satisfying the hypothesis of Lemma \ref{rescdiff}, we can apply the previous density estimates to $p_F$.  
Then the change of variable $y=\hat{x}_0+ \bar{A}_\d z$ gives that,  for $|y-\hat{x}_0|_{\bar{A}_\d(x_0)} \leq r$, we obtain respectively
\ben
\item
$\quad \frac{1}{ C |\det \bar{A}_\d(x_0)|}
\exp\left(-C |y-\hat{x}_0|_{\bar{A}_\d(x_0)}^2\right)\leq p_\d(x_0,y)$
\item $\quad p_\d(x_0,y)
\leq 
\frac{e^C }{ |\det\bar{A}_\d(x_0)|}\exp\left(-C^{-1} |y-\hat{x}_0|_{\bar{A}_\d(x_0)}^2\right)$
\een
where  $p_\d(x_0,y)$ is the density of $X_\d$ in $y$. These estimates and the equivalence between $|\cdot|_{A_\d}$ and  $|\cdot|_{\bar{A}_\d}$ (see \eqref{eq3} in the appendix) imply the thesis.
\epr

\br{e}
In the proof of Lemma \ref{rescdiff} 
we have used \textbf{A2}, the assumption of uniformly bounded derivatives, to say that $n_{F,G,U_r}(p)\leq e^C$ and $\|R_\d\|_{2,p}\leq e^C \d^2$, $C\in \cal{C}$. If we also ask that
\be{boundedcoefficients}
|\s(x)|+|b(x)|\leq \rho,\quad \forall x\in \R^2
\ee
we have that $n_{F,G,U_r}\leq \tilde{C}$ and  $\|R_\d\|_{2,p}\leq \tilde{C} \d^2$, $\tilde{C}\in \cal{C}$. This holds because, supposing the boundedness of the coefficients, we do not need anymore to use the Gronwall lemma to estimate the moments, but a direct computation is enough. These are standard estimates. In particular, in \eqref{deltabound} we have $1/\bar{C}$ instead of $\exp(-\bar{C})$.
As a consequence, if we also suppose \eqref{boundedcoefficients},   the lower bound in Lemma \ref{rescdiff} and Theorem \ref{mtstime} holds for $\d^*\in 1/\cal{C}$. In particular, taking $r^*=(L\vee C)^{-1/2}$ in Theorem \ref{mtstime}-(1) we can state that: $\exists r^*,\d^*\in 1/\cal{C}$, $C\in \cal{C}$ such that for $\d\leq \d^*$, for 
$|y-\hat{x}_0|_{A_\d(x_0)} \leq r^*$
\[
\frac{1}{ C \d^2}
\leq p_\d(x_0,y)
\]
On the other hand, in the upper bound we cannot get rid of the exponential dependence in the constant. Indeed, the estimate on $\Gamma_F(64)$ of Lemma \ref{n} is involved (the estimate on the ``non-degeneracy'' of the rescaled diffusion $F$). This has an exponential dependence on the parameters, even supposing \eqref{boundedcoefficients}, because it involves the moments of $Z_t$, the inverse of the flow of $X$, and in this estimate we always need to use Gronwall lemma.
Anyways, taking $r^*=(L)^{-1/2}$ in Theorem \ref{mtstime}-(2) we find that:
$\exists r^*,\d^*\in 1/\cal{C}$, $C\in \cal{C}$ such that for $\d\leq \d^*$, for 
$|y-\hat{x}_0|_{A_\d(x_0)} \leq r^*$
\[
p_\d(x_0,y)\leq \frac{e^C}{  \d^2}
\]
We put together those two inequalities in the following two-sided bound, which is the formulation that will be used to prove the tube estimate:\\
$\exists r^*,\d^*\in 1/\cal{C}$, $C\in \cal{C}$ such that for $\d\leq \d^*$, for 
$|y-\hat{x}_0|_{A_\d(x_0)} \leq r^*$
\be{denco}
\frac{1}{ C \d^2}
\leq
p_\d(x_0,y)
\leq \frac{e^C}{  \d^2}.
\ee
\er

\section{Tube estimates of the diffusion process}\label{diffusion}
As an application of Theorem \ref{mtstime} we prove the tube estimate. 
We suppose in this section  $\s,\,b\in \C^5(\R^2)$ and set, for $x\in \R^2$,
\[
n(x)=\sum_{k=0}^5 \sum_{|\a|= k} |\partial^\a_x b(x)| + |\partial^\alpha_x \s(x)|,
\quad
\quad \l(x)=\l_*(A(x)).
\]
We consider the diffusion \eqref{eqn} on $[0,T]$, and the skeleton path \eqref{control}: for $\phi\in L^2[0,T]$, let
\[
x_t(\phi)=x_0+\int_0^t \s(x_s(\phi)) \phi_s ds + \int_0^t b(x_s(\phi)) ds,\quad \mbox{ for } t\in [0,T].
\]
Recall \textbf{H1}, \textbf{H2}, \textbf{H3}, \textbf{H4}: 
\[
\l(y)\geq \l_t, \quad\quad
n(y)\leq n_t,\quad\quad
\partial_\s \s(y)=\kappa_\s(y) \s(y),\quad\quad \forall |y-x_t(\phi)|<1, \quad \forall t\in[0,T]
\]
Moreover, defining $(R_t)_{ t\in[0,T]}$ the time-dependent radius of the tube, we suppose that
\begin{align*}
n_\cdot &:[0,T]\rightarrow [1,\infty)& &R_\cdot :[0,T]\rightarrow (0,1]&& \\ 
\l_\cdot &:[0,T]\rightarrow (0,1] & &|\phi_\cdot|^2 :[0,T]\rightarrow (0,\infty)&&
\end{align*}
are in $\in L(\mu,h)$, for some $h>0,\, \mu\geq 1$, where $L(\mu,h)$ is the class of non-negative functions which have the property
\[
f(t)\leq \mu f(s) \quad\mbox{ for }  |t-s| \leq h.
\]
Denote, for $0\leq t\leq T$, for $K_*$, $q_*$ positive universal constants,
\be{Rmax}
R_t^*(\phi)=\exp\left(-K_* \left(\frac{\mu n_t}{\l_t}\right)^{q_*}\mu^{2q_*}  \right)
\left( h\wedge \inf_{0\leq \d\leq h} 
\left\{
\d \big/ \int_t^{t+\d} |\phi_s|^2 ds \right\} 
\right)
\ee
\bt{mttubes}
Let $(X_t)_{t\in [0,T]}$ be a process verifying \eqref{eqn}, and $x_t(\phi)$ the skeleton path defined above. If  {\bf H1}, {\bf H2}, {\bf H3}, {\bf H4} are satisfied,
there exist positive universal constants $\bar{K},\bar{q}$ such that
\[
\exp\left(- \int_0^T \bar{K} \left(\frac{\mu n_t}{\l_t}\right)^{\bar{q}}\left( \frac{1}{h}+\frac{1}{R_t}+|\phi_t|^2 dt\right) \right)\leq
\PR\left(\sup_{t\leq T} |X_t-x_t(\phi)|_{A_{R_t}(x_t(\phi))}\leq 1 \right).
\]
Moreover, there exist positive universal constants $\bar{K},\bar{q},K_*,q_*$  such that if $R_. \leq R_. ^*(\phi)$
\[
\begin{split}
&\PR\left(\sup_{t\leq T} |X_t- x_t(\phi)|_{A_{R_t}(x_t(\phi))}\leq 1 \right)\\
&\quad\quad\leq \exp\left(- \int_0^T \bar{K} \left(\frac{\mu n_t}{\l_t}\right)^{\bar{q}}
\left(\frac{\exp\left(-K_* \left(\frac{\mu n_t}{\l_t}\right)^{q_*}\right)}{R_t}+|\phi_t|^2\right)dt \right)
\end{split}
\]
\et
\br{ooo}
Remark that for $R_t \leq R_t ^*(\phi)\leq h\, \exp\left(-K_* \left(\frac{\mu n_t}{\l_t}\right)^{q_*}\right)$ the statement in \eqref{tuberesult} is implied by this one.
\er
\bpr
A main point in this proof is the choice a sequence of short time intervals in a way that allows us to apply the short time density estimate. This issue is related to the choice of a an ``elliptic evolution sequence" in \cite{BallyKohatsu:10,[BFM]}.
We fix $\phi$ from the beginning and write $x_t$ for $x_t(\phi)$ to have a more readable notation.

We introduce also the time-dependent version of  \eqref{def:constant}. For $t\in[0,T]$
\be{def:constant:time}
\cal{C}_t=\left\{C_t>0: C_t=\exp\left( K\left(n_t/\l_t)^q\right)\right),\exists K,q\geq 1 \right\}
\ee
The constants defined above depend on $\s,b$ through the ratio $n_t/\l_t$ locally along the skeleton path. We stress that $K,q$ do not depend on $\s,b$ and do not depend on $t\in [0,T]$.
We will also denote by $1/\cal{C}_t=\{\d_t>0:1/\d_t\in \cal{C}_t\}$. \\
We start proving the lower bound.
\\
STEP 1 (Time grid and notations):  
We set, for large $q_1,K_1$ to be fixed in the sequel, 
\[
f_R(t)= K_1 \left(\frac{\mu n_t}{\l_t}\right)^{q_1} \left( \frac{1}{h}+\frac{1}{R_t}+|\phi_t|^2\right).
\]
We use this function to split the time interval $[0,T]$ is short-enough sub-intervals (our time grid).
Recall \textbf{H4}:  $|\phi _{.}|^{2},n_{.},\l_{.}, R_{.}\in
L(\mu ,h)$, $\exists \mu\geq 1$, $0<h\leq 1$.
This implies $f_R\in L(\mu^{2 q_1+1},h)$. We also define
\be{tdeltaw}
\d(t)=\inf_{\d>0} \left\{ \int_t^{t+\d} f_R(s) ds \geq \frac{1}{\mu^{2 q_1+1}} \right\}.
\ee
Since 
\[
\frac{\d(t)}{h} = \int_t^{t+\d(t)} \frac{1}{h} ds \leq \int_t^{t+\d(t)} f_R(s)ds =\frac{1}{\mu^{2 q_1+1}},
\]
for any $t\in [0,T]$, $\d(t)\leq h/\mu^{2 q_1+1}\leq h$. Therefore we can use on the intervals $[t,t+\d(t)]$ the fact that our bounds are in $L(\mu,h)$. If $0<t-t'\leq h$,
\[
\mu^{2q_1+1} f_R(t) \d(t) \geq \int_{t}^{t+\delta(t)} f_R (s) ds=
\frac{1}{\mu^{2q_1+1}}
=\int_{t'}^{t'+\d(t')}f_{R}(s)ds \geq \mu^{-(2q_1+1)} f_{R}(t)\d(t'),
\]
so $\d(t')/\d(t)\leq \mu^{4q_1+2}$. Also the converse holds, and
$\d(\cdot) \in L(\mu^{4q_1+2},h)$. We set
\[
\ve(t)=\left( \int_t^{t+\delta(t)} |\phi_s|^2 ds\right)^{1/2}.
\]
We have
\[
\frac{1}{\mu^{2 q_1+1}}= \int_t^{t+\d(t)} f_R(s) ds 
\geq \int_t^{t+\d(t)} \frac{f_R(t)}{\mu^{2 q_1+1}} ds
\geq \d(t) \frac{f_R(t)}{\mu^{2 q_1+1}},
\]
so 
\be{Rdeltaw}
\d(t) \leq \frac{1}{f_R(t)} \leq \frac{R_t}{K_1} \left(\frac{\l_t}{\mu n_t}\right)^{q_1}.
\ee
Similarly,
\[
\frac{1}{\mu^{2 q_1+1}} 
\geq \int_t^{t+\d(t)} K_1 \left(\frac{\mu n_s}{\l_s}\right)^{q_1} |\phi_s|^2 ds 
\geq \frac{1}{\mu^{2 q_1}}  K_1 \left(\frac{\mu n_t}{\l_t}\right)^{q_1} \ve(t)^2, 
\]
and we can write both 
\be{energyw}
\d(t) \leq \frac{1}{K_1} \left(\frac{\l_t}{\mu n_t}\right)^{q_1},\quad\mbox{ and }\quad
\ve(t)^2 \leq \frac{1}{K_1} \left(\frac{\l_t}{\mu n_t}\right)^{q_1}.
\ee
We set our time grid as
\[
t_0=0;\quad t_k=t_{k-1}+\d(t_{k-1}),
\]
and introduce the following notation on the grid:
\[
\d_k=\d(t_k);\quad \ve_k=\ve(t_k);\quad n_k=n_{t_k};\quad \l_k= \l_{t_k};\quad X_k=X_{t_k};\quad x_k=x_{t_k};\quad R_k=R_{t_k}.
\]
We also define
\[
\hat{X}_k= X_k + b(X_k)\delta_k;
\quad
\hat{x}_k=x_k+b(x_k)\delta_k,
\]
and for $t_k\leq t \leq t_{k+1}$,
\[
\hat{X}_k(t)=X_k+b(X_k)(t-t_k);
\quad
\hat{x}_k(t)=x_k+b(x_k)(t-t_k).
\]
Moreover we denote 
\[
|\xi|_k=|\xi|_{A_{\delta_k}(x_k)};\quad \cal{C}_k=\cal{C}_{t_k},
\]
and $r^*_k\in \cal{C}_k$ the radius $r^*$ associated to \eqref{denco}, when taking as initial condition $x_0=x_k$.

\br{}
Consider $D_k=\{ \sup_{t_k\leq t \leq t_{k+1}} |X_t-x_t|_{A_{R_t}(x_t)}\leq 1 \},$ and
$\Gamma_k = \{ |X_k-x_k|_k \leq r_k \}$, where $r_k$ is radius smaller than $1$ that will be defined in the sequel. We denote $\PR_k$ the conditional probability
\[
\PR_k(\cdot)=\PR\left( \cdot|W_t,t\leq t_k;  \Gamma_k \right)
\]
We will lower bound
$\PR\left(\sup_{t\leq T} |X_t-x_t(\phi)|_{A_{R_t}(x_t(\phi))}\leq 1 \right)$
computing the product of the probabilities $\PR_k \left(D_k\cap \Gamma_{k+1}\right) $, and this computation relies on the application of the density estimate in short time. Remark that {\bf A1}, {\bf A3} are local assumption, therefore it is enough to ask for {\bf H1}, {\bf H3} to apply Theorem \ref{mtstime}. What about {\bf A2} (global) and {\bf H2} (local)? Suppose that we have a process $X$ which, for some external reasons, verifies \eqref{eqn} for $t_k\leq t \leq t_{k+1}$, and such that $\sup_{t_k\leq t \leq t_{k+1}} |X_t-x_t|_{A_{R_t}(x_t)}\leq 1 $. From ${\bf H2}$ 
\[
n(y)\leq n_k \quad\mbox{ for } \{y\in \R^2: |y-x_k|\leq 1\}
\]
A classical theorem (see \cite{whitney1944}) tells us that we can define $\bar{\s},\bar{b}$ which coincide with $\s,b$ on $\{y\in \R^2: |y-x_k|\leq1\}$, which are differentiable as many times as $\s,b$ but on the whole $\R^2$, and for which
\[
n(y)\leq \a n_k \quad\mbox{ for all }  y\in \R^2,  \mbox{ with } \a \mbox{ constant. }
\]
Let $\bar{X}$ be the strong solution to 
\[
\bar{X}_t=X_k+\int_{t_k}^t \bar{\s}(\bar{X}_s)\circ dW_s + \int_{t_k}^t \bar{b}(\bar{X}_s) ds, \quad t\in[t_k,t_{k+1}].
\]
It is clear that 
\[
\PR ( D_k  \cap \Gamma_{k+1})=
\PR \big(  \{ \sup_{t_k\leq t \leq t_{k+1}} |\bar{X}_t-x_t|_{A_{R_t}(x_t)}\leq 1 \} \cap \{ |\bar{X}_{t_{k+1}}-x_{k+1} |_{k+1} \leq r_{k+1} \}\big),
\]
and therefore we can equivalently prove our estimates supposing that $n(y)$ is globally, and not just locally, bounded by $n_k$. From now on we assume that
$n(y)\leq n_k$ for $y\in \R^2$.
\er

\noindent
STEP 2 (Application of the density estimate):
Lemmas \ref{normcomp}, \ref{normpoint}, \ref{normcontrolW}, \ref{normcontrolWf}  hold for $\d_k$ and $\ve_k$ small enough, and in particular Lemma \ref{normcontrolWf} says that
\be{ci0}
\frac{1}{C^1_k}|\xi|_{A_\d(x_k)} \leq |\xi|_{A_\d(x_{k+1})} \leq C^1_k |\xi|_{A_\d(x_k)},
\ee
for some $C^1_k\in\cal{C}_k$, for any $\d\leq \d_k$. 
Recall \eqref{energyw}, and 
\[
R_k/\mu\leq R_t\leq \mu R_k,\quad \mbox{ for } t_k\leq t \leq t_{k+1},
\]
so that $R_t\geq \d_k$ for $t_k\leq t \leq t_{k+1}$.
Moreover we have $|x_{k+1}-\hat{x}_k|_k\leq C_k(\ve_k\vee\sqrt{\d_k})$, and for all $t_k\leq t\leq t_{k+1}$, applying also \eqref{basicmatrnorm1}, $|x_t-\hat{x}_k(t)|_{A_{R_t}(x_t)} \leq C_k (\ve_k\vee\sqrt{\d_k})$ for $t_k\leq t \leq t_{k+1}$. Recall again \eqref{energyw}, and we fix $q_3,K_3$ such that, for $q_1\geq q_3,K_1\geq K_3$, the Lemmas \ref{normcomp}, \ref{normpoint}, \ref{normcontrolW}, \ref{normcontrolWf} hold and
\begin{align}
|x_{k+1}-\hat{x}_k|_k &\leq r^*_k/8 \label{magia1w}\\
|\hat{x}_k(t)-x_t|_{A_{R_t}(x_t)} &\leq \frac{1}{4} \mbox{ for all } t_k\leq t\leq t_{k+1} \label{magia2w},
\end{align}
and moreover $\d_k\leq \d^*_k$ associated to \eqref{denco} with initial condition $x_k$.

Now, $\d(\cdot) \in L(\mu^{4q_1+2},h)$ implies $\d_k/\d_{k+1}\leq 
\mu^{4q_1+2}$ and $\d_{k+1}/\d_k\leq \mu^{4q_1+2}$.  This, \eqref{ci0} and \eqref{basicmatrnorm1} give
\be{few}
\frac{1}{C^1_k \mu^{2q_1+1}}|\xi|_k \leq |\xi|_{k+1}\leq \mu^{2q_1+1}C^1_k |\xi|_k,
\ee
where $C^1_k$ is in $\cal{C}_k$, depending on $K_3,q_3$. We now set, for $K_2,q_2$ to be fixed in the sequel,
\be{rayw}
r_k= \frac{1}{K_2 \mu^{2q_1+2q_2+1} }\left(\frac{\l_k}{n_k}\right)^{q_2},
\ee
and define as we said before
\[
\Gamma_k=\{ |X_k-x_k|_k \leq r_k \},\quad D_k=\{ \sup_{t_k\leq t \leq t_{k+1}} |X_t-x_t|_{A_{R_t}(x_t)}\leq 1 \},
\]
and $\PR_k$ as the conditional probability
\[
\PR_k(\cdot)=\PR\left( \cdot|W_t,t\leq t_k; \Gamma_k \right).
\]
We find a lower bound for $\PR_k(\Gamma_{k+1}\cap D_k)$ using our density estimate in short time. 
We denote $p_k(X_k,y)=p_{\d_k}(X_k,y)$ the density of $X_{k+1}$ in $y$ with respect to $\PR_k$.
We prove that on $\{y:|y -x_{k+1}|_{k+1}\leq r_{k+1}\}$ we can apply
\eqref{denco} to $p_k(X_k,\cdot)$ and so there exists $\underline{C}_k\in \cal{C}_k$ such that
\be{eqstttw}
\frac{1}{\underline{C}_k \d_k^2} \leq 
p_{k}(X_k,y)
\ee
We estimate
\be{decow}
|y-\hat{X}_{k}|_k \leq
|y-x_{k+1}|_k+|x_{k+1}-\hat{x}_k|_k+|\hat{x}_k-\hat{X}_k|_k.
\ee
We already have \eqref{magia1w}. Since we are on $|y -x_{k+1}|_{k+1}\leq r_{k+1}$, from \eqref{few} and the fact that $r_{k+1}/r_k\leq \mu^{2 q_2}$
\[
|y-x_{k+1}|_k 
\leq C_k^1 \mu^{2 q_1+1} |y-x_{k+1}|_{k+1}
\leq C_k^1 \mu^{2 q_1+1} r_{k+1}
\leq C_k^1 \mu^{2 q_1+2 q_2 + 1} r_{k}
\leq  \frac{C_k^1}{K_2}\left(\frac{\l_k}{n_k}\right)^{q_2}. 
\]
It also holds $|\hat{x}_k-\hat{X}_k|_k \leq C_k |{x}_k-{X}_k|_k \leq C_k r_k$, for some $C_k\in\cal{C}_k$. Similarly, since $R_t\geq \d_k$, from \eqref{basicmatrnorm1}
$|\hat{x}_k(t)-\hat{X}_k(t)|_{A_{R_t}(x_t)} \leq C_k r_k$, for all $t_k \leq t \leq t_{k+1}$. Recalling \eqref{rayw}, we can fix $K_2,q_2$ such that 
$|y-x_{k+1}|_k\leq r^*_k/16$, $|\hat{x}_k-\hat{X}_k|_k \leq r^*_k/16$, and
\be{magia3w}
|\hat{X}_k(t)-\hat{x}_k(t)|_{A_{R_t}(x_t)} \leq 1/4,\quad \mbox{ for all } t_k\leq t\leq t_{k+1}.
\ee
From \eqref{decow}, \eqref{magia1w} this implies $|y-\hat{X}_{k}|_k \leq r^*_k/4$. 
We also have $|x_k-X_k|_k \leq r_k$, so we can also fix $K_2,q_2$ such that $r_k \leq \a$ in Lemma \ref{normpoint}. Therefore 
\[
\frac{1}{4}|\xi|_k\leq|\xi|_{A_{\d_k}(X_k)}\leq 4|\xi|_k.
\]
So $|y-\hat{X}_{k}|_{A_{\d_k}(X_k)} \leq r^*_k$ and \eqref{denco} holds (which means that \eqref{eqstttw} holds). Now, from Lemma \ref{normpoint} and \eqref{few}
\[
\begin{split}
\{|\cdot -x_{k+1}|_{A_{\d_k}(X_k)}\leq r_{k+1}/(4C^1_k\mu^{2 q_1+1}) \}
& \subset \{|\cdot -x_{k+1}|_{k}\leq r_{k+1}/(C^1_k\mu^{2 q_1+1} )\} \\
& \subset \{|\cdot -x_{k+1}|_{k+1}\leq r_{k+1}\}   ,
\end{split}
\]
and $ r_{k+1}/(4 C^1_k\mu^{2 q_1+1} )\geq r_k /(4 C^1_k \mu^{2 q_1+ 2q_2+1} )= \frac{1}{ 4 C^1_k K_2 \mu^{4 q_1+ 4 q_2+2} }\left(\frac{\l_k}{n_k}\right)^{q_2}$. So
\[
Leb(|\cdot -x_{k+1}|_{k+1}\leq r_{k+1}) \geq 
\d^2_k \det A(X_k)  
\left(\frac{1}{ 4 C^1_k K_2 \mu^{4 q_1+ 4 q_2+2} }\left(\frac{\l_k}{n_k}\right)^{q_2}\right)^2.
\]
Now, from {\bf H1}, $\det A(X_k)\geq \l_k$. So, from \eqref{eqstttw}, 
\[
\PR_k(\Gamma_{k+1})\geq \frac{1}{\underline{C}_k} \left(\frac{1}{ 4 C^1_k K_2 \mu^{4 q_1+ 4 q_2+2} }\left(\frac{\l_k}{n_k}\right)^{q_2}\right)^2 \l_k
\]
where $\underline{C}_k \in \cal{C}_k $ is the constant in \eqref{denco}. This implies
\[
2 \mu^{-8 q_1} \exp(-K_4 (\log \mu + \log n_k -\log \l_k )) \leq \PR_k(\Gamma_{k+1})
\]
for some constant $K_4$ (depending on $K_2,K_3,q_2,q_3$; on the contrary, we keep explicit the dependence in $q_1$, which is not fixed yet). 

\noindent
STEP 3 (Concatenation): Consider now $t_k\leq t \leq t_{k+1}$. Recall the definition
\[
D_k=\left\{ \sup_{t_k\leq t \leq t_{k+1}} |X_t-x_t|_{A_{R_t}(x_t)}\leq 1 \right\},
\]
and introduce 
\[
E_k=\left\{\sup_{t_k\leq t\leq t_{k+1}} |X_t-\hat{X}_k(t)|_{A_{R_t}(x_t)}\leq \frac{1}{2}\right\}.
\]
We decompose
\[
|X_t-x_t|_{A_{R_t}(x_t)}
\leq|X_t-\hat{X}_k(t)|_{A_{R_t}(x_t)}+|\hat{X}_k(t)-\hat{x}_k(t)|_{A_{R_t}(x_t)}
+|\hat{x}_k(t)-x_t|_{A_{R_t}(x_t)},
\]
and, from the previous part of the proof, \eqref{magia2w} gives $|\hat{x}_k(t)-x_t|_{A_{R_t}(x_t)} \leq 
1/4$, and \eqref{magia3w} gives $|\hat{X}_k(t)-\hat{x}_k(t)|_{A_{R_t}(x_t)} \leq 1/4$. So $|X_t-x_t|_{A_{R_t}(x_t)}\leq |X_t-\hat{X}_k(t)|_{A_{R_t}(x_t)}+1/2$, and therefore $E_k\subset D_k$. 

Now we have to estimate $E_k$.
A development of $X_t-\hat{X}_k(t)$ similar to \eqref{dec} gives that the diffusion moves with speed $\d_k^{1/2}$ in the direction of $\s(x_k)$, $\d_k^{3/2}$ otherwise. Taking the $|\cdot|_{A_{R_t}(x_t)}$ norm we account precisely of this fact. Applying the exponential martingale inequality we find that
\[
\PR_k(E_k^c) \leq  \exp\left(- \frac{1}{K_5} \left(\frac{\l_k}{\mu n_k}\right)^{q_5} \frac{R_k}{\d_k} \right)
\]
for some constants $K_5,q_5$. From \eqref{Rdeltaw}, $R_k/\d_k\geq K_1(\mu n_k/\l_k)^{q_1}$. We recall that $\l_k\leq 1$ and $n_k\geq 1$, so choosing and fixing now $q_1,K_1$ large enough we conclude
\[
\PR_k(E^c_k) \leq \mu^{-8 q_1} \exp(-K_4 (\log \mu + \log n_k -\log \l_k )) \leq \frac{1}{2}\PR_k(\Gamma_{k+1}),
\]
so
\be{shorttimediffw}
\begin{split}
\PR_k(\Gamma_{k+1}\cap D_k) &\geq \PR_k(\Gamma_{k+1}\cap E_k)\geq \PR_k(\Gamma_{k+1})-\PR_k(E^c_k)\geq \frac{1}{2} \PR_k(\Gamma_{k+1})\\
&\geq \exp\left( -K_6 (\log \mu +\log n_k -\log \l_k) \right),
\end{split}
\ee
for some constant $K_6$.
Let now $N(T)=\max\{k: t_k\leq T\}$. From Definition \eqref{tdeltaw}
\[
\int_0^T f_R(t)dt \geq \sum_{k=1}^{N(T)} \int_{t_{k-1}}^{t_k} f_R(t)dt \geq  
\frac{N(T)}{\mu^{2q_1+1}}.
\]
From \eqref{shorttimediffw},
\begin{align*}
\PR\left(\sup_{t\leq T} |X_t-x_t(\phi)|_{A_{R_t}(x_t(\phi))}\leq 1 \right)
& \geq \PR \left(\bigcap_{k=1}^{N(T)} \Gamma_{k+1}\cap D_k\right) \\
& \geq \prod_{k=1}^{N(T)} \exp(-K_6(\log \mu +\log n_k -\log \l_k))\\
& = \exp\left(-K_6\sum_{k=1}^{N(T)}\log \mu +\log n_k -\log \l_k\right).
\end{align*}
Since
\[
\begin{split}
\sum_{k=1}^{N(T)}(\log \mu +\log n_k -\log \l_k)
&=\mu^{2q_1+1}\sum_{k=1}^{N(T)}\int_{t_k}^{t_k+1} f_R(s)ds (\log \mu +\log n_k -\log \l_k)\\
&\leq \int_0^T \mu^{2q_1+1} f_R(t) \log\left(\frac{\mu^{3} n_t}{\l_t}\right)dt,
\end{split}
\]
the lower bound follows.

\noindent
STEP 4 (Upper bound): We define, with the same $K_1,q_1$ as in STEP 1,
\[
g_R(t)= 
K_1 \left(\frac{\mu n_t}{\l_t}\right)^{q_1}
\left( \frac{\exp \left( -K_*\left(\frac{\mu n_t}{\l_t}\right)^{q_*} \mu^{2q_*} \right)}{R_t} 
+|\phi_t|^2 \right)
\]
Because of \eqref{Rmax}, for all $t\in[0,T$],
\be{rh}
\frac{\exp \left( -K_*\left(\frac{\mu n_t}{\l_t}\right)^{q_*}\mu^{2q_*} \right)}{R_t} \geq \frac{1}{h}
\ee
We define now a new $\d(t)$
\[
\d(t)=\inf_{\d>0} \left\{ \int_t^{t+\d} g_R(s) ds \geq \frac{1}{\mu^{2q_1+1}} \right\}
\]
and, as before,
\[
\ve(t)=\left( \int_t^{t+\delta(t)} |\phi_s|^2 ds\right)^{1/2}.
\]
As in STEP 1, using also \eqref{rh}, we can check that \eqref{energyw} holds also for this choice of $\d$:
\[
\d(t) \leq \frac{h}{K_1} \left(\frac{\l_t}{\mu n_t}\right)^{q_1} \leq \frac{1}{K_1} \left(\frac{\l_t}{\mu n_t}\right)^{q_1},\quad\mbox{ and }\quad
\ve(t)^2 \leq \frac{1}{K_1} \left(\frac{\l_t}{\mu n_t}\right)^{q_1}.
\]
In particular, $\d(t)\leq h$. With these definitions we set a time grid $\{t_k: k=0,\dots, N(T)\}$ and all the associated quantities as in STEP 1.  As we did for the lower bound, since we estimate the probability of remaining in the tube for any $t\in [t_k,t_{k+1}]$, we can suppose that the bound $n(y)\leq n_k$ holds $\forall y\in \R^2$.
The short time density estimate \eqref{denco} holds again.
Recall now that $R_{.}\in L(\mu,h)$, and this gives the analogous to \eqref{few}:
\[
\frac{1}{C^1_k \sqrt{\mu}}|\xi|_{A_{R_k}(x_k)}\leq
|\xi|_{A_{R_{k+1}}(x_{k+1})} \leq
C^1_k \sqrt{\mu} |\xi|_{A_{R_k}(x_k)}
\]
We define
\[
\Delta_k=\{ |X_k-x_k|_{A_{R_k}(x_k)} \leq 1 \},
\]
$\tilde{\PR}_k$ as the conditional probability
 $\tilde{\PR}_k(\cdot)=\PR\left( \cdot|W_t,t\leq t_k;  \Delta_k \right)$.
Now, since $\d(t)\leq h$, we can apply the fact that $R,\l,n\in L(\mu,h)$ and
\[
\begin{split}
\int_t^{t+\d(t)}  K_1 \left(\frac{\mu n_s}{\l_s}\right)^{q_1} |\phi|_s^2 ds 
\leq
\mu^{2 q_1} K_1 &\left(\frac{\mu n_t}{\l_t}\right)^{q_1} \int_t^{t+\d(t)}   |\phi|_s^2 ds,  \\
\int_t^{t+\d(t)}  K_1 \left(\frac{\mu n_s}{\l_s}\right)^{q_1}
\frac{\exp \left( -K_*\left(\frac{\mu n_s}{\l_s}\right)^{q_*}\mu^{2q_*}\right) }{R_s} ds & \\
  \leq
\mu^{2q_1+1} K_1 &\left(\frac{\mu n_t}{\l_t}\right)^{q_1}
\exp \left( -K_*\left(\frac{\mu n_t}{\l_t}\right)^{q_*} \right) \frac{\d(t)}{R_t}.
\end{split}
\]
Recall now \eqref{Rmax}
\[
R_t\leq R_t^*(\phi)=\exp\left(-K_* \left(\frac{\mu n_t}{\l_t}\right)^{q_*} \mu^{2q_*}\right)
\left( \inf_{0\leq \d\leq h} 
\left\{
\d \big/ \int_t^{t+\d} |\phi_s|^2 ds \right\} 
\right),
\]
which implies
\[
\int_t^{t+\d(t)} |\phi_s|^2 ds\leq
\exp\left(-K_* \left(\frac{\mu n_t}{\l_t}\right)^{q_*} \right)\frac{\d(t)}{R_t} 
\]
We obtain
\[
1=
\mu^{2q_1+1}\int_t^{t+\d(t)}  g_R(s) ds \leq
2\mu^{4q_1+2}
K_1 \left(\frac{\mu n_t}{\l_t}\right)^{q_1}
\exp \left( -K_*\left(\frac{\mu n_t}{\l_t}\right)^{q_*} \right) \frac{\d(t)}{R_t}
\]
so
\be{Rdelta}
\frac{R_t}{\d(t)}
\leq
2\mu^{4q_1+2}
K_1 \left(\frac{\mu n_t}{\l_t}\right)^{q_1}
\exp \left( -K_*\left(\frac{\mu n_t}{\l_t}\right)^{q_*} \right)
\ee
As we did in STEP 1, if $q_*,K_*$ are large enough, $R_{k}$ is small enough and the upper bound for the density holds on $\Delta_{k+1}$. 
Because of \eqref{ci0},
\[
Leb( |\cdot-x_k|_{A_{R_k}(x_{k+1})} \leq 1) \leq Leb( |\cdot-x_k|_{A_{R_k}(x_k)}\leq 1) (C^1_k)^2 = (C^1_k)^2\det(A(x_k)){R_k}^2.
\]
Now, using the density estimate, 
\[
\tilde{\PR}_k(\Delta_{k+1})\leq
(C^1_k)^2\det(A(x_k))\,e^{\overline{C}_k} 
\left(\frac{{R_k}}{\d_k}\right)^2.
\]
where $\overline{C}_k$ is the constant in the upper bound of \eqref{denco}. 
Recall \eqref{Rdelta}, for $t=t_k$
\[
\frac{R_k}{\d_k}
\leq
2\mu^{4q_1+2}
K_1 \left(\frac{\mu n_k}{\l_k}\right)^{q_1}
\exp \left( -K_*\left(\frac{\mu n_k}{\l_k}\right)^{q_*} \right)
\]
so we chose now $K_*,q_*$ large enough to have
\[
\tilde{\PR}_k(\Delta_{k+1})\leq \exp(-K_7)
\]
for a constant $K_7>0$. From the definition of $N(T)$
\[
\int_0^T g_R(t)dt = \sum_{k=1}^{N(T)} \int_{t_{k-1}}^{t_k} g_R(t)dt = \frac{N(T)}{\mu^{2q_1+1}}\leq N(T).
\]
As before
\begin{align*}
& \PR\left(\sup_{t\leq T} |X_t-x_t(\phi)|_{A_{R_t}(x_t(\phi))}\leq 1 \right)\leq \prod_{k=1}^{N(T)} \tilde{\PR}_k(\Delta_{k+1})\\
& \quad \leq \prod_{k=1}^{N(T)} \exp(-K_7)= \exp(-K_7\,N(T))\leq 
\exp\left(-K_7\int_0^T g_R(t)\right),
\end{align*}
and we have the upper bound. 
\epr

\section{Matrix norm and control metric}

\subsection{Matrix norms}\label{matrnorm}
In this work we use a number of properties of norms associated to the matrix $A$ and $A_R$. Recall that in general we can associate a norm to a matrix $M$ with full row rank via
\[
|y|_M=\sqrt{\langle (MM^T)^{-1}y,y  \rangle}.
\]
Recall that, for $R>0$,
\[
A=\left(\s,[\s,b]\right),\quad \quad \quad A_R=\left(R^{1/2}\s,R^{3/2}[\s,b]\right)
\]
\bl{basicmatrnorm}
For every $y\in \mathbb{R}^2$ and $0<R\leq R' \leq 1$,
\begin{eqnarray}
& (R/R')^{1/2} |y|_{A_R}\geq |y|_{A_{R'}}\geq (R/R')^{3/2} |y|_{A_R} \label{basicmatrnorm1}\\
& \frac{1}{R^{1/2}\l^*(A)}|y|\leq |y|_{A_R} \label{2} \leq \frac{1}{R^{3/2}\l_*(A)}|y| \label{basicmatrnorm2}
\end{eqnarray}
\el
\begin{proof}
Writing explicitly the inequalities \eqref{basicmatrnorm1}, we easily see that they are verified if $0< R\leq R'\leq1$. 
Taking $R'=1$, we have
$$
R^{1/2} |y|_{A_R}\geq |y|_{A}\geq R^{3/2} |y|_{A_R}
$$
and so
$$
\frac{1}{R^{1/2}\l^*(A)}|y| \leq |y|_{A_R} \leq \frac{1}{R^{3/2}\l_*(A)}|y|
$$
\end{proof}

\br{remmatr}
Recall the following properties of matrices: 
\begin{align*}
\forall \xi,\quad
C\, |\xi|_{B}^2\geq |\xi|_{A}^2  
\Leftrightarrow 
C\, \left( B B^T \right)^{-1} \geq \left( A A^T \right)^{-1}
\Leftrightarrow 
B B^T \leq C\,A A^T
\end{align*}
and, denoting with $M_i$ the columns of $M$,
\[
\langle MM^T\xi,\xi\rangle=\sum_i\langle M_i,\xi\rangle^2,
\]
so that
\[
\l_*(M)^2=\inf_{|\xi|=1}\sum_{i} \langle M_i,\xi\rangle^2\quad\mbox{ and }\quad
\l^*(M)^2=\sup_{|\xi|=1}\sum_{i} \langle M_i,\xi\rangle^2
\]
Taking $M=A(x)=(\s(x), [\s,b](x))$ we have in particular that
\begin{equation}\label{eigA}
\l_*(A(x))^2|\xi|^2 \leq \langle \sigma(x),\xi \rangle^2 + \langle [\s,b](x),\xi \rangle^2\leq \l^*(A(x))^2|\xi|^2\quad \forall \xi\in \R^2
\end{equation}
\er
We prove now some equivalences between norms that will be needed especially in the concatenation along the tube. We state them for $t_k=t_0=0$ to lighten the notation.
Recall that $x_0$ is the initial condition of \eqref{eqn}, and that in the concatenation (Section \ref{diffusion}) we have
\ben
\item[{\bf H1}]
$\l_*(A(x))\geq \l_0,\quad \forall |x-x_0|<1$
\item[{\bf H2}] $n(x)\leq n_0, \quad \forall x\in \R^2$ (this is justified in STEP 1 of the proof)
\item[{\bf H3}] $\partial_\s \s(x)=\kappa_\s(x) \s(x) ,\quad \forall |x-x_0|<1$, $|\kappa_\s|\leq n_0,\,|\nabla\kappa_\s |\leq n_0$
\een
Moreover, we recall that $\l_0\leq 1$ and $n_0\geq 1$.
In \eqref{def:constant:time} we define a class of constants that in the case $t=0$ is 
\[
\cal{C}_0=\left\{C>0: C=\left( K\left( n_0/\l_0)^q\right)\right),\exists K,q\geq 1 \right\}
\]
\bl{normcomp} 
There exists $C\in \cal{C}_0, \d^*\in 1/\cal{C}_0$ such that for $\d\leq \d^*$, with $\hat{x}_0=x_0+b(x_0)\d$, for any $\xi\in\mathbb{R}^2$
\beq
\label{eq3}
&\frac{1}{C}|\xi|_{A_\d(x_0)} \leq |\xi|_{\bar{A}_\d(x_0)}\leq C|\xi|_{A_\d(x_0)}\\
\label{eq4}
&\frac{1}{C} |\xi|_{A_\d(x_0)} \leq |\xi|_{A_\d(\hat{x}_0)} \leq C |\xi|_{A_\d(x_0)}
\eeq
\el

\br{}
This lemma is used also in Section \ref{sshortime}, when $\cal{C}_t$ has not yet been defined. It is clear that in that case the constants must be taken in $\C$ defined in \eqref{def:constant}.
\er
\bpr
We take $M=A_\d(x_0)$ and $M=\bar{A}_\d(x_0)$ in Remark \ref{remmatr}. Recall that $\l_0\leq 1$ and $n_0\geq 1$ and notice that
\[
|\partial_b \s(x_0)|\leq 4 n_0^2 \leq \frac{4n_0^2}{\l_*(A(x_0))}\l_*(A(x_0))\leq C \l_*(A(x_0)), \mbox{ with } C\in \cal{C}_0
\]
so, from \eqref{eigA}
\[
\d^3 \langle \partial_b\s(x_0),\xi\rangle^2
\leq \d^3 C \l_*^2(A(x_0)) |\xi|^2 \leq C (\d\langle\s(x_0),\xi \rangle^2 + \d^3\langle [\s,b](x_0),\xi \rangle^2).
\]
We have
\begin{align*}
\d \langle  \s(x_0) &+ \d\partial_b\s(x_0),\xi \rangle^2 + \d^3\langle  [\s,b](x_0),\xi \rangle^2 \\
& \leq 2 \d \langle  \s(x_0),\xi \rangle^2 + 2 \d^3 \langle \partial_b\s(x_0),\xi \rangle^2 + 
\delta^3\langle [\s,b](x_0),\xi \rangle^2\\
& \leq C (\d\langle\s(x_0),\xi \rangle^2 +
\d^3\langle[\s,b](x_0),\xi \rangle^2),
\end{align*}
so $|\xi|^2_{A_{\d}(x_0)}\leq C|\xi|^2_{\bar{A}_\d(x_0)}$. Analogously, 
\begin{align*}
\d \langle  \s(x_0),\xi \rangle^2 + \d^3\langle [\s,b](x_0),\xi \rangle^2\leq C (\langle \d\s(x_0)+ \d\partial_b\s(x_0),\xi \rangle^2+\d^3\langle [\s,b](x_0),\xi \rangle^2),
\end{align*}
so
$|\xi|_{\bar{A}_\d(x_0)}^2\leq C |\xi|^2_{A_{\d}(x_0)}$. From
\[
|\s(\hat{x}_0)-\s(x_0)|=|\s(x_0+b(x_0)\d)-\s(x_0)|\leq\int_0^\d |\nabla \s(x_0+b(x_0)t) b(x_0)|dt \leq C \d,
\]
applying again Remark \ref{remmatr} as in the previous point, also \eqref{eq4} follows.
\epr
The following lemma establish the equivalence of matrix norms of this kind when the matrix is  taken in two points that are close in such matrix norms. 
\bl{normpoint}
Consider $x_0,x, y\in\R^2$, with $|x-x_0|< 1$. There exist $\a\in 1/\cal{C}_0$ such that if and $|x-y|_{A_\d(x)}\leq \a$,
\[
\frac{1}{4}|\xi|_{A_\delta(x)} \leq |\xi|_{A_\delta(y)}\leq 4 |\xi|_{A_\delta(x)},\quad\quad \forall \xi\in \R^2
\]
\el
\begin{proof}
Remark that \eqref{basicmatrnorm2} implies 
\[
|x-y|\leq \d^{1/2} C_1 |x-y|_{A_\d(x)}\leq \a C_1 \delta^{1/2}\leq \d^{1/2}
\]
for $\a\leq 1/C_1$. A Taylor development gives
\[
\s(x)-\s(y)=\nabla\s (x) (x-y)+\cal{O}(|x-y|^2),
\]
so
\begin{align*}
\langle \s(x),\xi\rangle^2
\leq 4 \langle \s(y),\xi\rangle^2+
4 \langle \nabla\s (x) (x-y),\xi\rangle^2 
+C_2 |x-y|^4 |\xi|^2.
\end{align*}
Since $A_\delta(x)$ is invertible, 
\[
\nabla \s(x) (x-y)=\nabla \s (x) A_\d(x) A^{-1}_\d(x) (x-y).
\]
From Cauchy-Schwartz inequality and $| A^{-1}_\d(x) (x-y)|\leq  \a$, 
\[
\begin{split}
|\langle \nabla\s(x) (x-y),\xi\rangle|
&=|\langle  A^{-1}_\d(x) (x-y), (\nabla\s(x) A_\d(x))^T \xi\rangle|\\
&\leq \a |(\nabla\s(x) A_\d(x))^T \xi|
\end{split}
\]
We are supposing {\bf H3}, so $\partial_\s \s=\kappa_\s \s$ holds in $x$, and
\begin{align*}
\nabla\s (x) A_\d(x)&= \nabla\sigma (x) (\delta^{1/2}\sigma(x),\delta^{3/2}[\sigma,b](x))\\
&=(\d^{1/2}\kappa_\s(x)\s(x), \d^{3/2}\partial_{[\s,b]}\s(x)),
\end{align*}
so
\[
|(\nabla\s(x) A_\d(x))^T \xi|^2=\d \kappa_\s^2(x) \langle\s(x),\xi\rangle^2+
\d^{3} \langle \partial_{[\s,b]}\s(x),\xi\rangle^2
\]
and therefore
\[
\begin{split}
\langle \nabla\s (x) (x-y),\xi\rangle^2
&\leq \a^2 (\d \kappa_\s^2(x) \langle\s(x),\xi\rangle^2+
\d^{3} \langle \partial_{[\s,b]}\s(x),\xi\rangle^2)\\
&\leq C_3 \a^2 \d \langle\s(x),\xi\rangle^2+
C_3 \a^2\d^{3} |\xi|^2
\end{split}
\]
Now,
\[
C_2 |x-y|^4 |\xi|^2\leq C_2 C_1^4 \a^4 \d^2 |\xi|^2  
\]
So
\[
\langle \s(x),\xi\rangle^2 \leq
 4 \langle \s(y),\xi\rangle^2+
4C_3 \a^2 \d \langle\s(x),\xi\rangle^2+
4C_3 \a^2\d^{3} |\xi|^2+
C_2 C_1^4 \a^4 \d^2 |\xi|^2  
\]
Taking $\a\leq \frac{1}{8C_3 C_2 C_1^2}$, this implies
\[
\langle \s(x),\xi\rangle^2 \leq
 8 \langle \s(y),\xi\rangle^2+ \a \d^2 |\xi|^2.
\]
In the direction $[\s,b]$ we have $[\s,b](x)-[\s,b](y)=\cal{O}(|x-y|)$
\begin{align*}
\langle [\s,b](x),\xi\rangle^2 \leq 2\langle [\s,b](y),\xi\rangle^2+ C_4 |x-y|^2 |\xi|^2
\leq 2\langle [\s,b](y),\xi\rangle^2+ C_4 C_1^2 \a^2 \d  |\xi|^2.
\end{align*}
We take now $\a\leq 1/(C_4 C_1^2)$, and we conclude that
\[
\d \langle \s(x),\xi\rangle^2 + \d^3 \langle [\s,b](x),\xi\rangle^2 \leq
 8 \d \langle \s(y),\xi\rangle^2+ 2\d^3 \langle [\s,b](y),\xi\rangle^2+ 2 \a \d^3 |\xi|^2.
\]
Using now \eqref{eigA} and {\bf H1}, 
\[
|\xi|^2\leq C_5 (\langle \s(y),\xi\rangle^2+ \langle [\s,b](y),\xi\rangle^2)
\]
So taking $\a\leq 4/C_5$  we have
\[
\d \langle \s(x),\xi\rangle^2 + \d^3 \langle [\s,b](x),\xi\rangle^2 \leq
 16 \d \langle \s(y),\xi\rangle^2+ 16\d^3 \langle [\s,b](y),\xi\rangle^2.
\]
From Remark \ref{remmatr}
we have $|\xi|_{A_\d(x)}\leq 4 |\xi|_{A_\d(y)}$. The converse inequality follows from an analogous reasoning. Remark that all the conditions we need on $\a$ are satisfied taking $\a \in 1/\cal{C}_0$ small enough, since $|x-x_0|< 1$ and {\bf H1}, {\bf H2}, {\bf H3}.

\end{proof}

We prove now that moving along a control $\phi\in L^2[0,T]$ for a small time, the trajectory remains close to the initial point in the $A_\d$-norm. 
Define, for fixed $\d$,
\[
\ve=\left(\int_0^\d |\phi_s|^2ds\right)^{1/2}.
\]
Recall that we have
\[
x_t(\phi)=x_0+\int_0^t \s(x_s(\phi)) \phi_s ds + \int_0^t b(x_s(\phi)) ds.
\]
\bl{normcontrolW}
There exist $\d_*,\ve_*\in 1/\cal{C}_0$, $C\in \C_0$ such that if $\d\leq \d_*,\,\ve\leq \ve_*$ 
\[
|x_\d(\phi)-(x_0+b(x_0)\d)|_{\bar{A}_\d(x_0)}\leq C (\ve\vee \d^{1/2}).
\]
\el 

\bpr
Via computations analogous to Decomposition \ref{dec} it is possible to write
\[
x_\d(\phi)-(x_0+b(x_0)\d)=\bar{A}_\d(x_0) (G_\phi+\tilde{R}_{\phi,\d})
\]
where
\[
G_\phi=\Th_\phi+\tilde{\eta}_\d(\Th_\phi),\quad 
\Th_\phi =\left( \begin{array}{l}
\d^{-1/2} \int_0^\d \phi_s ds\\
\d^{-3/2} \int_0^\d (\d-s)\phi_s ds
\end{array}\right)
\]
and
\[
| \tilde{R}_{\phi,\d} | \leq C ( \ve \vee \d^{1/2}).
\]
Remark that, by H\"older inequality,
\[
|\d^{-1/2} \int_0^\d \phi_s ds|\leq\ve,\quad |\d^{-3/2} \int_0^\d (\d-s)\phi_s ds|\leq\ve
\]
so $|\Th_\phi|\leq 2 \ve$ and by \eqref{dereta} $|\tilde{\eta}_\d(\Th_\phi)|\leq 4\ve^2$. Therefore $|G_\phi|\leq 4 \ve$ and
\[
|\bar{A}_\d(x_0)^{-1} (x_\d(\phi)-(x_0+b(x_0)\d))|=|G_\phi+\tilde{R}_{\phi,\d}|\leq C (\ve\vee\d^{1/2}).
\]
\epr

\bl{normcontrolWf}
There exist $\d_*,\ve_*\in 1/\cal{C}_0$, $C\in \cal{C}_0$ such that for $\d\leq \d_*,\,\ve\leq \ve_*$ 
\[
\frac{1}{C}|\xi|_{A_\d(x_0)}
\leq |\xi|_{A_\d(x_\d)}
\leq C |\xi|_{A_\d(x_0)}
\]
\el 
\bpr
Recall $\hat{x}_0=x_0+\d b(x_0)$.
Applying in this order \eqref{eq4}, \eqref{eq3}, Lemma \ref{normcontrolW} we obtain
\[
|x_\d-\hat{x}_0|_{A_\d(\hat{x}_0)}
\leq C|x_\d-(x_0+b(x_0)\d)|_{A_\d(x_0)}
\leq C|x_\d-(x_0+b(x_0)\d)|_{\bar{A}_\d(x_0)}\leq C (\ve\vee \d^{1/2}).
\]
Now, choosing $\d_*,\ve_*$ small enough, we can apply Lemma \ref{normpoint} to the points $x_\d,\hat{x}_0$, and    
\[
\frac{1}{4}|\xi|_{A_\d(\hat{x}_0)}
\leq |\xi|_{A_\d(x_\d)}
\leq 4 |\xi|_{A_\d(\hat{x}_0)}.
\]
Now again \eqref{eq4} concludes the proof.
\epr

\subsection{The control metric}\label{controlmetric}
Recall \eqref{eqn}, \eqref{defA}, \eqref{defAR}. 
In the spirit of \cite{KusuokaStroock:87}, we want to express our results is some control norm. 
Let
\[
\Omega=\{x\in \R^2\,:\,\l(x)=\l_*(A(x))>0\}
\]
A natural way to associate a quasi-distance to the matrix norm $|\cdot|_{A_R(\cdot)}$ used in this paper is to define
\[
d(x,y) < \sqrt{R} \Leftrightarrow |x-y|_{A_R(x)} < 1.
\]
(we take $\sqrt{R}$ because it is the``diffusive" regime). With this definition, $d$ is a quasi-distance on $\Omega$,  verifying the following properties (see \cite{NagelSteinWainger:85}):
\begin{itemize}
\item[$i)$] for every $x\in \Omega$, for every $r>0$, the set $\{y\in\Omega\,:\,d(x,y)<r\}$ is open
\item[$ii)$] $d(x,y)=0$ if and only if $x=y$
\item[$iii)$] for every compact set $K\Subset\Omega$ there exists $C>0$ such
that $d(x,y)\leq C\big(d(x,z)+d(z,y)\big)$ holds for every $x,y,z\in K$ 
\end{itemize}
We say that two quasi-distances $d_{1}:\Omega\times \Omega\rightarrow \R^{+}$ and $%
d_{2}:\Omega\times \Omega\rightarrow \R^{+}$ are equivalent if for every
compact set $K\Subset\Omega$ there exists a constant $C$ such that for every 
$x,y\in K$ 
\begin{equation}
\frac{1}{C}d_{1}(x,y)\leq d_{2}(x,y)\leq Cd_{1}(x,y).  \label{Norm7a}
\end{equation}%
In particular if $d_{1}$ is a distance and $d_{2}$ is equivalent with $d_{1}$ then $d_{2}$ is a quasi-distance. 

On the other hand, the distance usually considered in the framework of hypoelliptic stochastic differential equations is the \emph{control distance} defined as follows: denote, for $x,y\in \Omega$,
\be{ODEsigma}
C(x,y)=\{\phi\in L^2(0,1): 
dv_s= \s(v_s) \phi_s ds,\,x=v_0,\,y=v_1
\}.
\ee
The control distance $d_c$ between $x$ and $y$ is
\[
d_c(x,y)=\inf \left\{\left(\int_0^1|\phi_s|^2 ds\right)^{1/2}:\phi\in C(x,y)   \right\}.
\]
Geometrically speaking, this corresponds to take the geodesic (i.e. the length-minimizing curve) joining $x$ and $y$ on the sub-Riemannian manifold associated with the diffusion coefficient $\s$. In our case this notion looks inadequate: we are supposing just a weak H\"{o}rmander condition, and this means that we have to use the drift coefficient $b$ to generate the whole space $\R^2$. Therefore any reasonable associated distance should incorporate $b$ as well. Moreover it should account of the different speed associated to the vector field $[\s,b]$. This is the reason for the following definition.
\bd{def:cara}
We first introduce a function which accounts of the different scale of propagation in the direction $[\s,b]$. For $\phi=(\phi^1_s,\phi^2_s)\in L^2((0,1),\R^2)$,
\[
\|\phi\|_{1,3}^2=
\int_0^1 |\phi_s^1|^2 ds + \left( \int_0^1 |\phi_s^2|^2 ds \right)^{\frac{1}{3}}
\]
We generalize \eqref{ODEsigma} to
\[
C_A(x,y)=\{\phi\in L^2((0,1),\R^2): 
dv_s= A (v_s) \phi_s ds,\,x=v_0,\,y=v_1
\}.
\]
A classic result by Carath\'eodory says that for any $x,y\in \Omega$ there exist a piecewise smooth $\phi\in C_{A}(x,y)$. We set
\[
d_c(x,y)=\inf 
\left\{ 
\|\phi\|_{1,3}
:\phi \in C_A(x,y)\right\}
\]
\ed
We are interested in establishing an equivalence between $d$, the quasi-distance defined via the matrix-norm, and $d_c$, the quasi-distance in terms of the control. 
\bl{g}
Let $\xi\in \Omega$.
Suppose that there exists a neighborhood $U_\xi$ of $\xi$ such that for all $x\in U_\xi$:
\ben
\item[{\bf A1'}] $\l_*(A(x))>\underline{\l}_\xi>0$,
\item[{\bf A2'}]  $\sum_{0\leq |\a| \leq 5}
|\partial_x^{\a} \s(x)|+|\partial_x^{\a} b(x)|
\leq \r_\xi ,$
\item[{\bf A3'}]  $\partial_\s \s(x)=\kappa_\s(x) \s(x),$ where $\kappa_\s$ is a differentiable scalar function, $|\kappa_\s(x)|\leq \rho_\xi$ and $|\nabla \kappa_\s(x)|\leq \rho_\xi$.
\een
Then there exist a neighborhood $V_\xi$ of $\xi$ and a constant $C_\xi$ such that, for any $x,y\in V_\xi$
\be{}
\frac{1}{C_\xi} d(x,y)\leq d_c(x,y)\leq C_\xi d(x,y).
\ee{}%
\el
\br{43}
This implies, using the fact that every open cover of a compact has a finite subcover, Corollary \ref{cor:cartu}. Moreover, again via a standard compactness argument, we have that \emph{if {\bf A1'}, {\bf A2'}, {\bf A3'} hold for any $\xi\in \Omega,$ then $d$ and $d_c$ are equivalent quasi-distances} on $\Omega$.
\er
\bpr
We use in this proof some notions on similar metrics and pseudo-metrics for which we refer to \cite{NagelSteinWainger:85}. For any $\phi \in L^\infty((0,T),\R^2)$ we set
\[
\|\phi\|_{1,3,\infty}=
\sup_{0\leq s\leq 1} |\phi_s^1|  + \sup_{0\leq s\leq 1} |\phi_s^2|^{\frac{1}{3}}
\]
and define
\[
\r(x, y) = 
\inf 
\left\{ 
\|\phi\|_{1,3,\infty}:\phi \in C_A(x,y)\right\}
\]
It is also possible to allow only constant linear combinations of the vector fields:
\be{v}
\bar{C}_A(x,y)=\{\th\in \R^2: 
dv_s= A (v_s) \th ds,\,x=v_0,\,y=v_1
\}
\ee
Analogously, we define
\[
\r_2(x, y) = 
\inf 
\left\{ 
|\th^1|+|\th^2|^{1/3}
:\,\th \in \bar{C}_A(x,y)\right\}
\]
In \cite{NagelSteinWainger:85} the quasi-distances $\r$ and $\r_2$ are defined in a slightly different way, but clearly equivalent to ours. It is also proved that  $\r$ and $\r_2$ are locally equivalent. We use here only the trivial inequality $\r\leq \r_2$. 
Remark that the difference between $\r$ and $d_c$ is that we take $\|\phi\|_{1,3,\infty}$ instead of $\|\phi\|_{1,3}$, so $d_c\leq \rho$ follows easily from the fact that the $L^2(0,1)$ norm is dominated by the $L^\infty(0,1)$ norm. 

We prove that, for fixed $\xi$, there exist $V_\xi$ and $C_\xi$ such that
\[
d(x,y)< \sqrt{R} \Rightarrow \r_2(x,y)< C_\xi\sqrt{R},
\] 
for $x,y\in V_\xi$. Since $x,y\in V_\xi$, we can suppose $|x-y|< \g_\xi$ small. By definition, $d(x,y)< \sqrt{R}$ means $|x-y|_{A_R(x)}< 1$. We prove that this implies the existence of $\th\in \bar{C}_A(x,y)$ with $|\th^1|< C_\xi R^{1/2},\,|\th^2|< C_\xi R^{3/2}$. Indeed, for fixed $x$, consider the function
\[
\th \rightarrow
\Phi(\th)= \int_0^1 A(v_s) \th ds,
\]
with $v$ satisfying $dv_s= A (v_s) \th ds$, $v_0=x$. Remark that $\Phi: \R^2\rightarrow \R^2$, $\Phi(0)=0$ and $\nabla \Phi (0)=A(x)$, which is non-degenerate because of {\bf A1'}. Therefore it is locally invertible: there exist two neighborhoods of $0$ such that $\Phi$ is a diffeomorphism from one to the other, and therefore for $y-x$ in the neighborhood in the image we can find $\th$ such that $\Phi(\th)=y-x$. Moreover, from the fact that {\bf A1'} and {\bf A2'} are uniform around $\xi$, the size of the neighborhoods can be taken uniformly in $x$. Therefore we can find a neighborhood of $\xi$ such that for given $x,y$ in this neighborhood, there exist $\th$ for which $\Phi(\th)=y-x$. Again from {\bf A1'} and {\bf A2'}, we can also suppose that $|\th|\leq C^\xi_1|\Phi(\th)|$. So, there exists $V_\xi$ neighborhood of $\xi$ such that, for $x,y\in V_\xi$, there exists $\th\in \bar{C}_A(x,y)$, and moreover
\[
|\th|\leq C_1^\xi |\Phi(\th)|\leq C_1^\xi |x-y| < C_1^\xi \g_\xi.
\]
We now show 
\[
|\th^1|< C_\xi R^{1/2},\,|\th^2|< C_\xi R^{3/2}.
\]
It is clear from\eqref{basicmatrnorm2} that $|\th^1|\leq |\th|\leq C_1^\xi |x-y|<C^\xi_2 R^{1/2}$. Now, with a development similar to \eqref{devn}, we can write
\[
\int_0^1 A(v_s)\th ds=[\s,b](x) \th^2+\s(x)\th^1+\eta(\th^1)+ L(\th),
\]
with $|L(\th)|\leq C_3^\xi|\th| \,(|\th^2|+|\th^1|^3)$ for $|\th| < C_1^\xi \g_\xi $ and $\eta$ defined as in \eqref{defeta}:
\[
\eta(u)= \left( \frac{\kappa_\s(x) }{2} u^2 + 
\frac{(\partial_\s \kappa_\s+\kappa_\s^2)(x)}{6} u^3 \right)\s(x)=
\left( \a(x)u^2+\b(x) u^3 \right)\s(x)
\]
(we have used again {\bf A3'}).
So
\[
A(x)^{-1} \int_0^1 A(v_s)\th ds=
\left(
\begin{array}{c}
\th^1+\a(x)(\th^1)^2+\b(x)(\th^1)^3\\
\th^2
\end{array}
\right)
+
A(x)^{-1} L(\th)
\]
Since $|\th|< C_1^\xi \g_\xi$ and $|\th^1|\leq C_2^\xi R^{1/2}$,
\[
|A(x)^{-1} L(\th)|\leq C_4^\xi|\th| \, (|\th^2|+|\th^1|^3 )\leq C_4^\xi C_1^\xi (C_2^\xi)^3\g_\xi (|\th^2|+R^{3/2} )\leq \frac{ |\th^2|+ R^{3/2}}{2},
\]
 choosing $\g_\xi\leq (2 C_4^\xi C_1^\xi(C_2^\xi)^3)^{-1}$.
In particular, the second component of $A(x)^{-1} L(\th)$ is in absolute value smaller 
than $(|\th^2|+R^{3/2} )/2$. Then the second component of $A(x)^{-1}\int_0^1 A(v_s)\th ds$ is in absolute value larger than $|\th^2|-(|\th^2|+R^{3/2} )/2 =(|\th^2|-R^{3/2} )/2$. As a consequence, the second component of $A_R(x)^{-1}\int_0^1 A(v_s)\th ds$ is in absolute value larger than $R^{-3/2}(|\th^2|-R^{3/2} )/2$. 
Since $|\int_0^1 A(v_s)\th ds|_{A_R(x)}=|x-y|_{A_R(x)}\leq 1$, we have $R^{-3/2}(|\th^2|-R^{3/2} )/2\leq 1$ and so we conclude $|\th^2|\leq 3R^{3/2}$.

We now prove
\[
d_c(x,y)< \frac{\sqrt{R}}{C^\xi_6} \Rightarrow d(x,y)< \sqrt{R}.
\]
We suppose $\phi\in C_A(x,y)$ with $\|\phi\|_{1,3}\leq \frac{\sqrt{R}}{C^\xi_6}$, which implies 
\[
\int_0^1 |\phi^1_s|^2 ds\leq \frac{ R 
}{(C^\xi_6)^2}
\quad \mbox{ and } \quad
\int_0^1 |\phi^2_s|^2  ds \leq\frac{ R^3 }{(C^\xi_6)^6}
\]
Developing as before and applying {\bf A3'},
\[
|x-y|_{A_{R}(x)} = 
\left|A_{R} (x)^{-1} \int_0^1 A(v_s)\phi_s  ds \right| \leq
C^\xi_7 \sqrt{ \frac{\int_0^1 |\phi^1_s|^2 ds}{R}+\frac{\int_0^1 |\phi^2_s|^2  ds
+(\int_0^1 |\phi^1_s|^2 ds )^3
}{R^3}}
\]
Therefore
\[
|x-y|_{A_{R}(x)} \leq 
C^\xi_7 \frac{ \sqrt{3}}
{C_6^\xi}
<1
\]
if $C^\xi_6$ is a large enough constant.
\epr

\section*{Acknowledgements}  
I thank Vlad Bally for introducing me to this subject and for his constant advice. I also thank Noufel Frikha, Arturo Kohatsu-Higa and St\'ephane Menozzi for stimulating discussions. I am grateful to the referees for their suggestions, which helped to improve the quality of the paper.

\bibliographystyle{plain} 
\bibliography{bibliografia}

\end{document}